\documentclass[a4paper,12pt]{article}
\usepackage{amsmath,amsthm,amscd,graphics,ifthen}
\newcommand{\ZZ}       {\ensuremath{\mathbf{Z}}}  
\newcommand{\CC}       {\ensuremath{\mathbf{C}}}  
\newcommand{\PP}       {\mathbf{P}}               
\newcommand{\CP}       {\mathbf{CP}}              
\newcommand{\p}        {\prime}   
\newcommand{\nbd}      {
   \nobreakdash-\hspace{0pt}}
\newcommand{\I}        {{[\,0,1\,]}}   
\newcommand{\set}[1]   {{\{\,#1\,\}}}  
\newcommand{\class}[1] {{[\,#1\,]}}    
\newcommand{\id}       {\ensuremath{\mathrm{id}}}   
\newcommand{\incl}     {\ensuremath{\mathrm{incl}}} 
\newcommand{\can}      {\ensuremath{\mathrm{can}}}  
\newcommand{\sing}     {\mathrm{sing}} 

\newcommand{\HPP}      {\mathbf{\hat P}} 
\newcommand{\wn}       {\mathrm{wn}}  
\newcommand{\pr}       {\mathrm{pr}}  
\newcommand{\Ll}       {\mathcal{L}}  
\newcommand{\Ee}       {\mathcal{E}}  
\newcommand{\Mm}       {\mathcal{M}}  
\newcommand{\Dd}       {\mathcal{D}}  
\newcommand{\Vv}       {\mathcal{V}}  
\newcommand{\shrpp}    {{\sharp\,\p}} 
\newcommand{\lstar}   {\,*}      
\newcommand{\lbullet} {\,\bullet}
\newcommand{\lhash}   {\,\#}     

\newcommand{\harrow}[5][]{
  \begin{picture}(0,0)\put(#2,#3){\raisebox{0pt}[0pt][0pt]%
   {\makebox[0pt]{$\xrightarrow[#1]{\makebox[#4pt]{$\scriptstyle #5$}}$}}}%
  \end{picture}}

\newcommand{\darrow}[8][]{
  \begin{picture}(0,0)\put(#2,#3)%
   {\rotatebox{#4}{%
   \raisebox{0pt}[0pt][0pt]{\makebox[0pt]{$\xrightarrow{\makebox[#5pt]{%
   \ifthenelse{\equal{#1}{is}}{$\scriptstyle\sim$}{}}}$}}}%
   \begin{picture}(0,0)\put(#6,#7){$\scriptstyle#8$}\end{picture}}%
  \end{picture}}

\newcommand{\isarrow}[1][]{
   \overset{#1}{\overset{\sim}{\longrightarrow}}}

\newcommand{\leftisarrow}[1][]{
   \overset{#1}{\overset{\sim}{\longleftarrow}}}
     
\newcommand{\his}[1]{
   \underset{\scriptstyle\sim}{#1}}

\newcommand{\shis}[1]{
   \smash[t]{\underset{\scriptstyle\sim}{#1}}}

\newcommand{\vis}{\wr\;}   

\DeclareMathOperator{\im}{im}   
\DeclareMathOperator{\coker}{coker}   
\newtheorem{prop}{Proposition}[section]
\newtheorem{theo}[prop]{Theorem}
\newtheorem*{antheo}{Theorem}

\newtheorem{coro}[prop]{Corollary} 
\newtheorem{lem}[prop]{Lemma} 

\theoremstyle{definition}
\newtheorem{dfn}[prop]{Definition}   

\newtheorem{notn}[prop]{Notation}

\theoremstyle{remark}
\newtheorem{rem}[prop]{Remark} 
\newtheorem{claim}[prop]{Claim}

\numberwithin{equation}{section}  

\begin{document}

\date{February 23, 2002}  

\title{Zariski-van~Kampen theorem for higher homotopy groups.}

\author{D. Ch\'eniot\\
\small Laboratoire d'Analyse, Topologie et Probabilit\'es (UMR CNRS 6632)\\
\small Centre de Math\'ematiques et Informatique, Universit\'e de Provence\\
\small 39, rue F. Joliot-Curie, 13453 Marseille Cedex 13, France\\
\small e-mail: cheniot@gyptis.univ-mrs.fr
\and A. Libgober\\
\small Department of Mathematics, University of Illinois at Chicago\\
\small 851 S.Morgan, Chicago, Illinois 60607, USA\\
\small e-mail: libgober@math.uic.edu}

\maketitle

\begin{abstract}
This paper gives an extension of the classical Zariski-van~Kampen
theorem describing the fundamental groups of the complements of plane
singular curves by generators and relations. It provides a procedure
for computation of the first non-trivial higher homotopy groups of
the complements of singular projective hypersurfaces in terms of
the homotopy variation operators introduced here.

\end{abstract}

\section{Introduction}
The classical Zariski-van~Kampen theorem expresses the fundamental
group of the complement of a plane algebraic curve in~$\CP^2$ as a
quotient of the fundamental group of the intersection of this
complement and a generic element of a pencil of lines
(cf.~\cite{Zariski}, \cite{vanKampen} and~\cite{DenisvanKampen}).
The latter group is always free and the quotient is taken by the
normal closure of a set of elements described in terms of the
monodromy action arising as a result of moving the above generic
element around the special elements of the pencil. This theorem is
the main tool for the study of the fundamental groups of the
complements of plane algebraic curves (cf.~\cite{survey}).

The purpose of the present paper is to describe a high dimensional
generalization of this theorem. Let $V$~be a hypersurface
of~$\CP^{n+1}$ having degree~$d$ and the dimension of its singular
locus equal to~$k$. It is shown in~\cite{Annals} that, if
$n-k\geq2$, $\pi_1(\CP^{n+1}-V)=\ZZ/d\ZZ$ and
$\pi_i(\CP^{n+1}-V)=0$ for $2\leq i\leq n-k-1$. Moreover the group
$\pi_{n-k}(\CP^{n+1}-V)$ depends on the local type and the position
of the singularities of~$V$. The latter homotopy group is called
the \emph{first non trivial (higher) homotopy group} of the
complement to~$V$ in~$\CP^{n+1}$. Since by the Zariski-Lefschetz
hyperplane section theorem (cf.~\cite{ZaLe}), for a generic linear
subspace~$H$ of codimension~$k$ in~$\CP^{n+1}$ one has
$\pi_{n-k}(\CP^{n+1}-V)=\pi_{n-k}(H-H\cap V)$, it is enough to
consider only the key case when $V$~has only isolated
singularities. This remark reduces also the case $n-k=1$ to the
Zariski-van~Kampen theorem mentioned above.

An analogue of the Zariski-van~Kampen theorem for higher homotopy
groups of the complement to hypersurfaces with isolated
singularities in $\CC^{n+1}$ was given in~\cite{Annals}. There it
was shown that, for a generic hyperplane~$L$, the homotopy group
$\pi_n(\CC^{n+1}-V)$ is the quotient of $\pi_n(L-L\cap V)$ by a
$\pi_1(L-L\cap V)$-submodule which depends not just on the
monodromy around the singular members of the pencil containing the
hyperplane section but also on certain ``degeneration operators" on
the homotopy groups of the special members of the pencil.

The present work proposes a different approach to a high
dimensional Zariski-van Kampen theorem. It is based on the
systematic use of homotopy variation operators introduced below.
Homological variation operators were considered in~\cite{London}
for a generalization of the second Lefschetz theorem
(cf.~\cite[Chap.~V, \S8, Th\'eor\`eme~VI]{Situs}, \cite{Wallace}
and~\cite{A-F}). From this point of view the main result of this
paper can also be viewed as a homotopy second Lefschetz theorem.

The main result of the paper is the following (restated as
Theorem~\ref{t:highvK} below):

\begin{antheo}
Let $V$~be a hypersurface in~$\PP^{n+1}$ with $n\geq2$ having only
isolated singularities. Consider a pencil $(L_t)_{t \in\PP^1}$ of
hyperplanes in~$\PP^{n+1}$ with  the base locus~$\Mm$ transversal
to~$V$. Denote by $t_1$, \dots,~$t_N$ the collection of those~$t$
for which $L_t\cap V$ has singularities. Let $t_0$~be different
from either of $t_1$, \dots,~$t_N$. Let $\gamma_i$~be a good
collection (cf.~Definition~\ref{t:loops}) of paths in~$\PP^1$
based in~$t_0$. Let $e\in\Mm-\Mm\cap V$ be a base point. Let
$\Vv_i$~be the variation operator (cf.~section~\ref{s:var})
corresponding to~$\gamma_i$. Then inclusion induces an isomorphism:
\begin{equation*}
 \pi_n(\PP^{n+1}-V,e)\leftisarrow\pi_n(L_{t_0}-L_{t_0}\cap V,e)\Big/
 \sum_{i=1}^N\im\Vv_i.
\end{equation*}
\end{antheo}

For $n=1$, this statement reduces to the classical
Zariski-van~Kampen theorem as we explain it in
Remark~\ref{generalize} below. But our proof does not work in the
case $n=1$. Thus we shall suppose $n\geq2$ in this article.

The paper is organized as follows. We start, in
section~\ref{s:prelim}, by describing in detail several pencils of
hyperplanes associated with a hypersurface~$V$ with isolated
singularities. We also review some vanishing results and
homological description of the homotopy groups of the complements
of hypersurfaces from~\cite{Annals}. In sections \ref{s:monodr},
\ref{s:degener} and~\ref{s:var} we describe the monodromies, the
degeneration operators and homotopy variation operators.
Section~\ref{s:link} describes the crucial relationship between
homotopy variation and degeneration operators. The last section
contains the announced theorem (theorem~\ref{t:highvK}). Two proofs
are given, one deriving it from the quoted theorem of~\cite{London}
and another from that of~\cite{Annals}.

Most of the work for this paper was done during several visits by
the first author to the University of Provence, to which he
expresses his gratitude for warm hospitality.  Both authors would
like to thank L\^e D\~ung Tr\'ang, David Trotman and Bernard
Teissier for reading parts of the manuscript and making valuable
comments.

Here is some notation we shall use throughout the paper.

\begin{notn}
\label{t:notn1}
\begin{enumerate}
\item The ground field in this paper is always~$\CC$ so we shall
omit~`$\CC$' from our notation of complex projective space which
becomes~$\PP^n$ for the
$n$\nbd dimensional one.

\item All inclusion maps will be denoted by~``$\incl$" and any
canonical surjection from a set to a quotient of it by~``$\can$". In
diagrams, we shall use the same letter for a map and any other map
obtained from it by restriction of the source or the target.
\label{t:notn1.maps}

\item All homology groups will be singular homology groups with integer
coefficients. Given a continuous map $f\colon X\to Y$ between topological
spaces, we shall denote by~$f_*$ the induced homomorphism
$H_n(X)\to H_n(Y)$, whatever be the integer~$n$, and by~$f_\bullet$ the
induced homomorphisms $C_n(X)\to C_n(Y)$ between chain groups. If
$X'\subset X$ and $f(X')\subset Y'\subset Y$, we shall write~$\bar f_*$ for
the induced homomorphisms $H_n(X,X')\to H_n(Y,Y')$. If $x\in X$ and
$y=f(x)\in Y$, we shall denote by~$f_\#$ the induced maps
$\pi_n(X,x)\to\pi_n(Y,y)$ for~$n\geq0$ and, if $x\in X'$ and
$y=f(x)\in Y'$, by~$\bar f_\#$ the induced maps
$\pi_n(X,X',x)\to\pi_n(Y,Y',y)$ for~$n\geq1$. All boundary operators in
homology or homotopy will be designated by~$\partial$. All absolute Hurewicz
homomorphisms will be denoted by~$\chi$ and the relative ones
by~$\bar\chi$. \label{t:notn1.topalg}

\item Given an absolute cycle~$\xi$, we shall write~$\class{\xi}_X$ for its
homology class in a space~$X$ containing it and, if $\eta$~is a chain
contained in~$X$ with boundary contained in $X'\subset X$, we shall denote
by~$\class{\eta}_{(X,X')}$ its homology class in~$X$ modulo~$X'$. If
$\xi'$~is an(other) absolute cycle contained in~$X$, we shall write
`$\xi\sim\xi'\text{ in X}$' to mean that $\xi$~is homologous to~$\xi'$
in~$X$. The homotopy class of a loop~$\gamma$ will be denoted
by~$\bar\gamma$, the space in which this class has to be taken being made
clear by the context. \label{t:notn1.Hom}

\item The singular locus of the algebraic hypersurface~$V$ will be
designated by~$V_\sing$.

\item In a blow up, the total transform of a subset~$E$ of the blown up
space will be denoted by~$\hat E$ and its strict transform by~$E^\sharp$.
\label{t:notn1.blow}
\end{enumerate}
\end{notn}

\section{Preliminaries}
\label{s:prelim}

\subsection{General setup}
\label{s:setting}

Let $V$~be a closed algebraic hypersurface of $(n+1)$\nbd
dimensional complex projective space~$\PP^{n+1}$, with only
\emph{isolated singularities}. Let $d$~be its degree. We suppose
$n\geq2$ for the reasons explained in the introduction and we
suppose $d\geq2$, the case $d=1$ being trivial (and appearing as a
combursome exceptional case in what follows). Let $M$~be a
projective $(n-1)$\nbd plane transverse to~$V$ (that is avoiding the
singular points of~$V$ and transverse to the non-singular part
of~$V$). Let $L$~be the pencil of hyperplanes with base locus~$M$,
that is, the set of projective hyperplanes of~$\PP^{n+1}$
containing~$M$. We want to compute the homotopy groups
of~$\PP^{n+1}-V$ with the help of its sections by the elements
of~$L$.

We take homogeneous coordinates $(x_1:\dots:x_{n+1}:x_{n+2})$
on~$\PP^{n+1}$, so chosen that $M$~is defined by the equations
\begin{equation*}
 x_{n+1}=x_{n+2}=0.
\end{equation*}
We then have a one-to-one parametrization of the elements of pencil~$L$ by
the elements of~$\PP^1$ as follows. Given also homogeneous coordinates
on~$\PP^1$, for each $t=(\lambda:\mu)\in\PP^1$, the hyperplane~$L_t$
of~$\PP^{n+1}$ with parameter~$t$ is defined by the equation
\begin{equation*}
 \lambda x_{n+1}+\mu x_{n+2}=0.
\end{equation*}
We shall thus consider~$L$ as being the parametrized
family~$(L_t)_{t\in{\PP^1}}$. The transversality of~$M$ to~$V$ entails the
following claim.

\begin{claim}
\label{t:transv}
The given choice of the axis~$M$ of pencil $L=(L_t)_{t\in{\PP^1}}$ is
generic. All the members of this pencil are transverse to~$V$ except a
finite number of them, say $L_{t_1}$, \dots,~$L_{t_N}$. Each
$L_{t_i}$~is transverse to~$V$ except at a \emph{finite} number of points,
which may be singular points of~$V$ or tangency points of~$L_{t_i}$ to the
non-singular part of~$V$, and moreover none of them belongs to~$M$.
\end{claim}

\begin{proof}
This is a consequence of~\cite[Corollaire~10.18, Corollaire~10.19
combined with Proposition~10.20 and Corollaire~10.17]{L'Ens}. The quoted
results apply when stratifying~$V$ by its singular
part~$V_\sing$ and non-singular part $V-V_\sing$. This is indeed a Whitney
stratification of~$V$ by Lemma~19.3 of~\cite{Wh}, since $V_\sing$~is $0$\nbd
dimensional.
\end{proof}
Thus, pencil~$(L_t)_{t\in{\PP^1}}$ looks like a stratified version
of the ``Lefschetz pencils'' of~\cite{SGA} but each $L_{t_i}\cap V$ may have
more than one singularity and these singularities may be of any kind.

Our goal is to define variation and degeneration operators on homotopy.
Each of those depends on a choice of~$L_{t_i}$ for some fixed index~$i$ 
and a loop~$\gamma_i$ running once around~$t_i$ in the parameter
space~$\PP^1$ and surrounding none of the points $t_1$, \dots,~$t_N$
besides~$t_i$. The main technical tool is an interpretation of the relevant
homotopy groups as the homology groups of universal coverings  which was
used in~\cite{Annals}. This material is discussed in the last part of this
section (cf.~\ref{s:homHom} below). The universal covers are obtained as
restrictions of a ramified cover of~$\PP^{n+1}$ by a hypersurface~$W$
of~$\PP^{n+2}$ viewed in the next subsection. The rest describes the 
classical blowing up construction in our framework which we use to get rid
of the base points of the pencils as will be needed for the definition of
degeneration operators.

\subsection{A ramified cover of $\PP^{n+1}$}
\label{s:ramcov}

In the homogeneous coordinates of~$\PP^{n+1}$ chosen in
section~\ref{s:setting}, let
\begin{equation*}
 f(x_1,\dots,x_{n+1},x_{n+2})=0
\end{equation*}
be an equation of~$V$ where $f$~is a homogeneous reduced polynomial of
degree~$d$.

Now, in~$\PP^{n+2}$ with homogeneous coordinates
$(x_0:x_1:\dots:x_{n+1}:x_{n+2})$, let $A$~be the point of coordinates
$(1:0:\dots:0:0)$. Let us consider the projection with center~$A$
\begin{equation*}
\begin{split}
 \pr\colon\quad \PP^{n+2}-\set{A}
  &\longrightarrow\PP^{n+1}\\
 (x_0:x_1:\dots:x_{n+1}:x_{n+2})
  &\longmapsto(x_1:\dots:x_{n+1}:x_{n+2}).
\end{split}
\end{equation*}
Let $W$~be the hypersurface of~$\PP^{n+2}$ given by the equation
\begin{equation*}
 x_0^d+f(x_1,\dots,x_{n+2})=0.
\end{equation*}
We have $A\notin W$. Thus
$\pi=\pr_{|W}$ is well defined. The following is a classical result.

\begin{claim}
\label{t:ramcov}
The map $\pi\colon W\to\PP^{n+1}$ is a holomorphic $d$\nbd fold covering 
of~$\PP^{n+1}$ totally ramified along~$V$.
\end{claim}

We consider also the embedding
\begin{equation*}
 j\colon\quad \PP^{n+1}\longrightarrow\PP^{n+2}\qquad
  (x_1:\dots:x_{n+2})\longmapsto(0:x_1:\dots:x_{n+2}),
\end{equation*}
the image of which is the projective hyperplane
$j(\PP^{n+1})\subset\PP^{n+2}$ given by $x_0=0$. We have
\begin{equation}
\label{jV}
 W\cap j(\PP^{n+1})=j(V)=\pi^{-1}(V),
\end{equation}
each of these subsets of~$\PP^{n+2}$ having equations equivalent to
\begin{equation*}
x_0=f(x_1,\dots,x_{n+2})=0.
\end{equation*}
The following claim is also easy to check from the equations.

\begin{claim}
\label{t:Wsing}
The singular points of~$W$ are the images by~$j$ of the singular points
of~$V$.
\end{claim}

Hypersurface~$W$ supports a natural pencil, the elements of which are
the  branched covers of the elements of pencil~$L$ and which can  be
explicitly described as follows. Let $\Mm$~be the projective $n$\nbd plane
of~$\PP^{n+2}$ defined by the same equations as~$M$ in~$\PP^{n+1}$, that is
\begin{equation*}
 x_{n+1}=x_{n+2}=0
\end{equation*}
and let
$\Ll=(\Ll_t)_{t\in{\PP^1}}$~be the pencil of
hyperplanes of~$\PP^{n+2}$ with base locus~$\Mm$. Here, with the
same homogeneous coordinates on~$\PP^1$ as in section~\ref{s:setting}, for
each $t=(\lambda:\mu)\in\PP^1$, the hyperplane~$\Ll_t$
of~$\PP^{n+2}$ is defined by the same equations as~$L_t$
in~$\PP^{n+1}$, that is
\begin{equation*}
 \lambda x_{n+1}+\mu x_{n+2}=0.
\end{equation*}
As a consequence we have
\begin{equation}
\label{piML}
 \pi^{-1}(M)=\Mm\cap W\qquad\text{and}\qquad
 \pi^{-1}(L_t)=\Ll_t\cap W\quad
 \text{for any $t\in\PP^1$}
\end{equation}
and also
\begin{equation*}
 \Mm\cap j(\PP^{n+1})=j(M)\quad\text{and}\quad
 \Ll_t\cap j(\PP^{n+1})=j(L_t)
 \text{\ for any $t\in\PP^1$}.
\end{equation*}

Unramified covers of $\PP^{n+1}-V$ and of its
sections by~$L_t$ and~$M$ are given by Claim~\ref{t:ramcov}. They can be
specified thanks to \ref{jV} and~\ref{piML} as follows.

\begin{claim}
\label{t:cov}
Map~$\pi$ induces the following holomorphic unramified $d$\nbd fold
coverings:
\begin{align*}
 \text{(i)}  &  &  W-j(V)            &\longrightarrow\PP^{n+1}-V,\\
 \text{(ii)} &  &  \Ll_t\cap(W-j(V)) &\longrightarrow L_t-L_t\cap V
    \quad\text{for any $t\in\PP^1$},\\
 \text{(iii)}&  &  \Mm\cap(W-j(V))   &\longrightarrow M-M\cap V.
\end{align*}
\end{claim}

It is worth noticing that pencil~$\Ll$ is good with respect to~$W$
and~$j(V)$ as $L$~was good with respect to~$V$. To make this precise, we
first stratify~$W$. The following claim is proved using
again~\cite[Lemma~19.3]{Wh}.

\begin{claim}
\label{t:Strat}
The partition
\begin{equation*}
 \Sigma=\set{j(V_\sing),j(V-V_\sing),W-j(V)}
\end{equation*}
of~$W$ is an algebraic Whitney stratification.
\end{claim}

The statement analogous to Claim~\ref{t:transv} is then the following. It
is a consequence of Claim~\ref{t:transv} that can be checked on the
equations.

\begin{claim}
\label{t:Transv}
The base locus~$\Mm$ of~$\Ll$ is transverse to the strata of~$\Sigma$ and
so is~$\Ll_t$ for all $t\in\PP^1$ distinct from $t_1$, \dots,~$t_N$.
Each~$\Ll_{t_i}$ is transverse to~$W-j(V)$, non-transverse to~$j(V_\sing)$
wherever it meets this finite set and transverse to~$j(V-V_\sing)$ except
at the points~$j(x)$ corresponding to the finite number of points~$x$ where
$L_{t_i}$~is tangent to the non-singular part of~$V$.
\end{claim}

\subsection{Blowing up the cover}
\label{s:Blow} 

The homotopical degeneration and variation operators we want will be
obtained by considering homological counterparts on the cover dealt with in
the preceding subsection. But the definition of the homological degeneration
operators will in turn require to blow up this cover along the base locus of
the pencil we considered. In fact we do it first for the ambient
space~$\PP^{n+2}$. Let
\begin{equation*}
 \HPP^{n+2}=\set{(x,\tau)\in\PP^{n+2}\times\PP^1\mid x\in \Ll_\tau}
\end{equation*}
be the blow up of~$\PP^{n+2}$ along~$\Mm$. It is given by the equation
\begin{equation*}
 \tau_1x_{n+1}+\tau_2x_{n+2}=0
\end{equation*}
which is separately homogeneous in the homogeneous coordinates
$(x_0:x_1:\dots:x_{n+1}:x_{n+2})$ of~$x$ and $(\tau_1:\tau_2)$ of~$\tau$.
This is an $(n+1)$\nbd dimensional complex analytic compact connected
submanifold of
$\PP^{n+2}\times\PP^1$.

The restrictions to~$\HPP^{n+2}$ of the projections
of~$\PP^{n+2}\times\PP^1$ onto its first and second factors give
respectively the blowing down morphism~$\Phi$ and projection~$P$:
\begin{equation*}
 \PP^{n+2}\overset{\Phi}{\longleftarrow}\HPP^{n+2}
 \overset{P}{\longrightarrow}\PP^1.
\end{equation*}
These are holomorphic mappings and $P$~is submersive.

For any subset $\Ee\subset\PP^{n+2}$, we shall, following
Notation~\ref{t:notn1}, item~\ref{t:notn1.blow}, denote
by~$\hat \Ee$ its total transform, that is
\begin{equation*}
 \hat \Ee=\Phi^{-1}(\Ee).
\end{equation*}
If $\Ee\subset \Mm$, then its total transform has a product
structure:
\begin{equation}
\label{Prod}
 \hat \Ee=\Ee\times\PP^1\quad
 \text{for any $\Ee\subset \Mm$}
\end{equation}
and the restrictions of $\Phi$ and~$P$ to~$\hat \Ee$ coincide with
the first and second projections. In particular, $\hat \Mm=\Mm\times\PP^1$.
The blowing down morphism induces an analytic isomorphism from
$\HPP^{n+2}-\hat \Mm$ onto $\PP^{n+2}-\Mm$. We shall be interested in the
total transform~$\hat W$ of the cover~$W$ of the preceding subsection.

For each~$t\in\PP^1$, we consider the strict transform of~$\Ll_t$ which we
denote by~$\Ll_t^\sharp$ following Notation~\ref{t:notn1},
item~\ref{t:notn1.blow}. We have
\begin{equation}
\label{Strict}
 \Ll_t^\sharp=\Ll_t\times\set{t}=P^{-1}(t).
\end{equation}
The blowing down morphism induces an analytic isomorphism
from~$\Ll_t^\sharp$ onto~$\Ll_t$ and hence an analytic isomorphism
\begin{equation}
\label{Isom}
 \Ll_t^\sharp\cap\hat \Ee\isarrow \Ll_t\cap \Ee\quad
 \text{for any $\Ee\subset\PP^{n+2}$}.
\end{equation}

At this point it will be convenient to introduce some abbreviations to save
space in diagrams.

\begin{notn}
\label{t:notn2}
We shall designate by~$E'$ the intersection of a subspace~$E$
of~$\PP^{n+1}$ (resp.~$\PP^{n+2}$, $\HPP^{n+2}$) with $\PP^{n+1}-V$
(resp.~$W-j(V)$, $\hat W-\widehat{j(V)}$). For instance,
$L_t'=L_t-L_t\cap V$, $\Mm'=\Mm\cap(W-j(V))$,
$\Ll_t^\shrpp=\Ll_t^\sharp\cap(\hat W-\widehat{j(V)})$ and even $\hat
W'=\hat W-\widehat{j(V)}$. We shall denote by~$P'$ the restriction of
projection~$P$ to
$\hat W-\widehat{j(V)}$.
\end{notn}

\subsection{Homological description of homotopy groups}
\label{s:homHom}

We obtain relevant homotopy groups of the base spaces of the covers
of Claim~\ref{t:cov} as homology groups of their total spaces except in
some exceptional cases where we shall content ourselves with a morphism from
a subgroup of the fundamental group of the base space onto the first
homology group of the total space.

\begin{lem}
\label{t:homHom}
We use Notation~\ref{t:notn2}. Let $e\in M'$
and~$\varepsilon\in\pi^{-1}(e)\subset\Mm'$ be base points.
\begin{enumerate}
\renewcommand{\theenumi}{(\roman{enumi})}
\renewcommand{\labelenumi}{\theenumi}
\item If~$n\geq2$, there are isomorphisms $\eta$ and~$\alpha_t$, for
 $t\in\PP^1-\set{t_1,\dots,t_N}$, defined by composition as follows:
 \begin{align*}
  \eta      &\colon\quad\pi_n(\PP^{n+1}-V,e)
     \isarrow[\pi_\#^{\smash[t]{-1}}]\pi_n(W-j(V),\varepsilon)
     \isarrow[\chi]H_n(W-j(V)),\\
  \alpha_t  &\colon\quad\pi_n(L_t',e)
     \isarrow[\pi_\#^{-1}]\pi_n(\Ll_t',\varepsilon)
     \isarrow[\chi]H_n(\Ll_t'),
 \end{align*}
 where the arrows labeled~$\chi$ are Hurewicz isomorphisms and the arrows
 labeled~$\pi_\#^{-1}$ are the inverses of isomorphisms induced by the
 projections of the coverings of Claim~\ref{t:cov} (which all are
 restrictions of map~$\pi$).
 \label{t:homHom.eta}
 Furthermore, for any $t\in\PP^1-\set{t_1,\dots,t_N}$, the following diagram
 is commutative:
 \begin{equation*}
 \begin{CD}
  H_n(\Ll_t')        @>\incl_*>>     H_n(W-j(V))\\
  @AA\vis\alpha_t A                  @AA\vis\eta A\\
  \pi_n(L_t',e)      @>\incl_\#>>    \pi_n(\PP^{n+1}-V,e)
 \end{CD}
 \end{equation*}
 (following Notation~\ref{t:notn1},
 item~\ref{t:notn1.maps}, all inclusion maps are denoted by~$\incl$).
\item If~$n\geq3$, there are isomorphisms~$\beta_i$, for $1\leq i\leq N$,
 and~$\gamma$ defined by composition as follows:
 \begin{align*}
  \beta_i  &\colon\quad\pi_{n-1}(L_{t_i}',e)
     \isarrow[\pi_\#^{\smash[t]{-1}}]\pi_{n-1}(\Ll_{t_i}',\varepsilon)
     \isarrow[\chi]H_{n-1}(\Ll_{t_i}'),\\
  \gamma   &\colon\quad\pi_{n-1}(M',e)
     \isarrow[\pi_\#^{-1}]\pi_{n-1}(\Mm',\varepsilon)
     \isarrow[\chi]H_{n-1}(\Mm'),
 \end{align*}
 where isomorphisms $\chi$ and~$\pi_\#^{-1}$ are as in item~(i).
 \label{t:homHom.beta3}
 Furthermore, for $1\leq i\leq N$, the following diagram is commutative: 
 \begin{equation}
 \label{gammabeta}
 \begin{CD}
  H_{n-1}(\Mm')  @>\incl_*>>
                                  H_{n-1}(\Ll_{t_i}')\\
      @AA\vis\gamma A                        @AA\vis\beta_i A\\
  \pi_{n-1}(M',e)              @>\incl_\#>>
                                  \pi_{n-1}(L_{t_i}',e).
 \end{CD}
 \end{equation}
\item When $n=2$, the projections of the coverings considered in
 item~\ref{t:homHom.beta3} induce isomorphisms
 \begin{align*}
  \pi_1(\Ll_{t_i}',\varepsilon)   &\isarrow[\pi_\#] G_i
                                      \subset\pi_1(L_{t_i}',e),\\
  \pi_1(\Mm',\varepsilon)         &\isarrow[\pi_\#] H
                                      \subset\pi_1(M',e),\\
 \end{align*}
 where $G_i$~is, for $1\leq i\leq N$, a subgroup of index~$d$ of
 $\pi_1(L_{t_i}',e)$ and $H$~a subgroup of index~$d$ of $\pi_1(M',e)$
 (the latter is free of rank~$d-1$). One can then define
 homomorphisms~$\beta_i$, for $1\leq i\leq N$,
 and~$\gamma$ by composition as follows:
 \begin{align*}
  \beta_i  &\colon\quad G_i
     \isarrow[\pi_\#^{\smash[t]{-1}}]\pi_1(\Ll_{t_i}',\varepsilon)
     \overset{\chi}{\longrightarrow}H_1(\Ll_{t_i}'),\\
  \gamma   &\colon\quad H
     \isarrow[\pi_\#^{\smash[t]{-1}}]\pi_1(\Mm',\varepsilon)
     \overset{\chi}{\longrightarrow}H_1(\Mm'),\\
 \end{align*}
 where the Hurewicz homomorphisms~$\chi$ are here abelianizations.
 Homomorphisms~$\beta_i$, for $1\leq i\leq N$, and~$\gamma$ are thus onto.
 \label{t:homHom.beta2}
 Furthermore, for $1\leq i\leq N$, the image of~$H$ by the natural map
 $\incl_\#\colon\pi_1(M',e)\to\pi_1(L_{t_i}',e)$ is included in~$G_i$ and
 the following diagram is commutative:
 \begin{equation}
 \label{gammabeta2}
 \begin{CD}
  H_1(\Mm')  @>\incl_*>>          H_1(\Ll_{t_i}')\\
      @AA\gamma A                        @AA\beta_i A\\
  H          @>\incl_\#>>         G_i.
 \end{CD}
 \end{equation}
\end{enumerate}
\end{lem}

\begin{proof}
Let $E\to B$ be one of the unramified coverings of Claim~\ref{t:cov}. Its
projection is a restriction of map~$\pi$. It induces an isomorphism from
the fundamental group of~$E$ onto a subgroup of index~$d$ of
the fundamental group of~$B$ and isomorphisms of
$k$\nbd th homotopy groups for $k\geq2$. But these vanish for
$2\leq k\leq n-1$ (this range may be empty) if $B=\PP^{n+1}-V$ because
$V$~is a hypersurface with isolated singularities (see
\cite[Lemma~1.5]{Annals}). The same is true if $B=L_t'=L_t-L_t\cap V$
with $t\neq t_i$ for $1\leq i\leq N$ because, by Claim~\ref{t:transv},
$L_t\cap V$ is a non-singular hypersurface of $L_t\simeq\PP^n$
(\cite[Lemma~1.1]{Annals}). On the other hand, the $k$\nbd th homotopy
groups vanish for $2\leq k\leq n-2$ (a range which may be empty) if
$B=L_{t_i}'$ with $1\leq i\leq N$ because then $L_{t_i}\cap V$ is, still by
Claim~\ref{t:transv}, a hypersurface with isolated singularities of
$L_{t_i}\simeq\PP^n$. And the same is true if $B=M'$ because $M$~was chosen
transverse to~$V$ so that $M\cap V$ is a non-singular hypersurface of
$M\simeq\PP^{n-1}$.

If moreover $E$~is simply connected, that is to say the cover universal,
then the Hurewicz homomorphism $\chi\colon\pi_n(E,\varepsilon)\to H_n(E)$
in the first two cases and
$\chi\colon\pi_{n-1}(E,\varepsilon)\to H_{n-1}(E)$ in the last two
will be an isomorphism due to the Hurewicz isomorphism theorem. By the next
lemma, this indeed happens for the values of~$n$ listed in the first two
items of the statement. In the cases of the last item, $\chi$~is an
epimorphism of abelianization. Besides when $n=2$, $M'$~is a projective
line minus $d$~points and its fundamental group is free of rank~$d-1$. In
the same case, the image of~$H$ by the natural map from $\pi_1(M',e)$ to
$\pi_1(L_{t_i}',e)$ is contained in~$G_i$ because the projections of the
coverings commute with inclusions. As to the commutativity of the diagrams,
it results from this and the functoriality of the Hurewicz homomorphisms.
\end{proof}

\begin{lem}
\label{t:univcov}
\begin{enumerate}
\renewcommand{\labelenumi}{(\roman{enumi})}
 \item If $n\geq2$, the first covering of Claim~\ref{t:cov} is a universal
  covering and so is the second one if $t\neq t_i$
  for~\hbox{$1\leq i \leq N$}.
 \item If~$n\geq3$, the third covering of Claim~\ref{t:cov} is a universal
  covering and so is the second covering with $t=t_i$ for~$1\leq i \leq N$.
\end{enumerate}
\end{lem}

\begin{proof}
Let $E\to B$ be one of these coverings.  According to the nature of
base space~$B$ as discussed in the proof of Lemma~\ref{t:homHom}, it is
pathwise connected and its fundamental group is~$\ZZ/d\ZZ$ when $n\geq2$ in
the first two cases there considered and when $n\geq3$ in the last two
(notice that all involved hypersurfaces have degree~$d$). Thus, for the
appropriate range of~$n$, this group has the same number of elements as the
fiber. The lemma then follows once it is verified that the total space~$E$
is pathwise connected. This can be seen from the irreducibility of~$W$,
$\Ll_t\cap W$ and~$\Mm\cap W$ or from the fact that, above a neighbourhood
of a regular point of~$V$, each of the coverings has a local model which is
a product of the cover associated with $z\mapsto z^d$ and a disk of
appropriate dimension.
\end{proof}

We shall also have to consider relative homotopy groups for the variation
operators.

\begin{lem}
\label{t:relhomHom}
Let $e$ and~$\varepsilon$ be base points as in Lemma~\ref{t:homHom}.
If~$n\geq 2$ and $t\in\PP^1-\set{t_1,\dots,t_N}$, there are
homomorphisms~$\bar\alpha_t$ defined by composition as follows (we use
Notation~\ref{t:notn2}):
\begin{equation*}
\bar\alpha_t\colon\quad\pi_n(L_t',M',e)
     \isarrow[\bar\pi_\#^{\smash[t]{-1}}]\pi_n(\Ll_t',\Mm',\varepsilon)
     \overset{\bar\chi}{\longrightarrow}H_n(\Ll_t',\Mm'),
\end{equation*}
where $\bar\chi$~is the relative Hurewicz homomorphism and
$\bar\pi_\#^{-1}$~the inverse of an isomorphism induced by the projection of
the second covering of Claim~\ref{t:cov}. Homomorphisms $\bar\chi$
and~$\bar\alpha_t$ are epimorphisms if~$n=2$ and isomorphisms if~$n\geq3$.
Furthermore, for~$n\geq3$ and $t\in\PP^1-\set{t_1,\dots,t_N}$, the following
diagram, where $\gamma$~is the isomorphism defined in
Lemma~\ref{t:homHom}, item~\ref{t:homHom.beta3} is commutative:
\begin{equation}
\label{alphagamma}
\begin{CD}
 H_n(\Ll_t',\Mm')  @>\partial>>  H_{n-1}(\Mm')\\
     @AA\vis\bar\alpha_t A              @AA\vis\gamma A\\
 \pi_n(L_t',M',e)  @>\partial>>  \pi_{n-1}(M',e).
\end{CD}
\end{equation}
When $n=2$, the image of the boundary homomorphism from $\pi_2(L_t',M',e)$
to $\pi_1(M',e)$ is contained in the subgroup~$H$ defined in
Lemma~\ref{t:homHom}, item~\ref{t:homHom.beta2} and the following
diagram, where $\gamma$~is the homomorphism defined there, is commutative:
\begin{equation}
\label{alphagamma2}
\begin{CD}
 H_2(\Ll_t',\Mm')  @>\partial>>  H_1(\Mm')\\
     @AA\bar\alpha_t A              @AA\gamma A\\
 \pi_2(L_t',M',e)  @>\partial>>  H.
\end{CD}
\end{equation}
\end{lem}

\begin{proof}
Map~$\pi$ induces the projection of the covering $\Ll_t'\to L_t'$ and we
have $\pi^{-1}(M')=\Mm'$ by \ref{jV} and~\ref{piML}. That
$\bar\pi_\#$~is then an isomorphism is a general fact about pairs of
fibered spaces (cf.~\cite[7.2.8]{Spanier}). We now come to the relative
Hurewicz homomorphism~$\bar\chi$. If~$n\geq2$ and
$t\in\PP^1-\set{t_1,\dots,t_N}$, spaces $\Mm'$ and~$\Ll_t'$ are
path-connected as seen in the proof of Lemma~\ref{t:univcov}.
Furthermore, the same lemma and the vanishing of higher homotopy groups
stated in the proof of Lemma~\ref{t:homHom} give then that
$\pi_k(\Mm',\varepsilon)=0$ for $0\leq k\leq n-2$ and
$\pi_k(\Ll_t',\varepsilon)=0$ for $1\leq k\leq n-1$, hence
$\pi_k(\Ll_t',\Mm',\varepsilon)=0$ for $1\leq k\leq n-1$ by the homotopy
exact sequence. The relative Hurewicz theorem (cf.~\cite[7.5.4]{Spanier})
then applies and gives that the Hurewicz homomorphism~$\bar\chi$ induces an
isomorphism onto $H_n(\Ll_t',\Mm')$ from the quotient of
$\pi_n(\Ll_t',\Mm',\varepsilon)$ obtained by identifying each element with
its images by the action of $\pi_1(\Mm',\varepsilon)$. But, by
Lemma~\ref{t:univcov}, this fundamental group is trivial if~$n\geq3$. When
$n=2$, the image of $\partial\colon\pi_2(L_t',M',e)\to\pi_1(M',e)$ is
contained in~$H$ because boundary homomorphisms commute with those induced
by $\pi$ and $\bar\pi_\#$~is onto as we have seen. As furthermore boundary
operators commute with Hurewicz homomorphisms (cf.~\cite[7.4.3]{Spanier}),
the last two diagrams are commutative.
\end{proof}

Finally the following lemma will be useful while computing homology in the
universal coverings.

\begin{lem}
\label{t:vanHom}
We have the following vanishing of homology groups.
\begin{enumerate}
\renewcommand{\labelenumi}{(\roman{enumi})}
\item $H_k(W-j(V))=0$ for $1\leq k\leq n-1$.
\item $H_k\bigl(\Ll_t\cap(W-j(V))\bigr)=0$ for $1\leq k\leq n-1$ if
 $t\neq t_i$ for~$1\leq i \leq N$.
\item $H_k\bigl(\Ll_{t_i}\cap(W-j(V))\bigr)=0$ for $1\leq k\leq n-2$
 and~$1\leq i \leq N$.
\item $H_k\bigl(\Mm\cap(W-j(V))\bigr)=0$ for $1\leq k\leq n-2$.
\end{enumerate}
\end{lem}

\begin{proof}
This also results from Lemma \ref{t:univcov}, the vanishing of higher
homotopy groups stated in the proof of Lemma~\ref{t:homHom} and the
Hurewicz isomorphism theorem.
\end{proof}
Note that the last two assertions are empty if~$n=2$.

\section{Monodromies}
\label{s:monodr}

The homological degeneration and variation operators which we must define at
the universal covering level involve monodromies around the exceptional
hyperplanes in the cover $W-j(V)$ and its blow up $\hat W-\widehat{j(V)}$.
These in turn are linked with a fibration structure of
$\hat W-\widehat{j(V)}$ outside of the exceptional hyperplanes, a structure
we shall also directly use for the degeneration operator.

\begin{claim}
\label{t:Fibr}
The restriction of~$P$
to $(\hat W-\widehat{j(V)})-\bigcup_{i=1}^N\Ll_{t_i}^\sharp$ is the
projection of a fiber bundle over
$\PP^1-\set{t_1,\dots,t_N}$. This bundle has
$\Mm\cap(W-j(V))\times(\PP^1-\set{t_1,\dots,t_N})$ as
a trivial subbundle of it. The fibers over
$t\in\PP^1-\set{t_1,\dots,t_N}$ are
$\Ll_t^\sharp\cap(\hat W-\widehat{j(V)})$ and
$\Mm\cap(W-j(V))\times\set{t}$.
\end{claim}

\begin{proof}
This results from the fact that, by Claim~\ref{t:Transv}, $\Mm$~is
transverse to the strata of a Whitney stratification of~$W$ for which
$W-j(V)$ is a union of strata: see~\cite[Corollary~3.12]{London}.
The quoted Corollary rests on the Thom-Mather first isotopy theorem.
\end{proof}
Notice however that there is not a similar fibration for $W-j(V)$ because
its sections by pencil~$(L_t)_{t\in\PP^1}$ are pinched together along
axis~$M$. Nevertheless we shall obtain monodromies also there, using the
isomorphisms from $\Ll_t^\sharp\cap(\hat W-\widehat{j(V)})$ to
$\Ll_t\cap(W-j(V))$ that, by~\eqref{Isom}, the blowing down morphism~$\Phi$
induces.

We shall consider the monodromies above a  special set of loops in the
parameter space $\PP^1$. We choose a base point~$t_0$ in~$\PP^1$ distinct
from points $t_1$, \dots,~$t_N$ and consider a \emph{good} system of
generators $(\bar\Gamma_i)_{1\leq i\leq N}$ of the fundamental group
$\pi_1(\PP^1-\set{t_1,\dots,t_N},t_0)$ (cf.~\cite{survey}
and~\cite[Definition~2.1]{Annals}). Recall that in such a system each
loop~$\Gamma_i$ is based at~$t_0$ and is the boundary of a subset~$D_i$
of~$\PP^1$ homeomorphic to a disk with $t_i$~as an interior point. Moreover
$D_i\cap D_j=\set{t_0}$ for~$i\neq j$. If the~$D_i$ are suitably chosen, a
presentation of $\pi_1(\PP^1-\set{t_1,\dots,t_N},t_0)$ is given by these
generators and the relation
$\bar\Gamma_1\dots\bar\Gamma_N=1$. A standard method to construct such a
system is to obtain it as the homotopy classes $\bar\gamma_i=\bar\Gamma_i$
of parametrized loops $\gamma_i\colon\I\to\PP^1-\set{t_1,\dots,t_N}$ as
described in the following definition.

\begin{dfn}
\label{t:loops}
Let $t_0$~be a base point in $\PP^1-\set{t_1,\dots,t_N}$. Let $\Delta_1$,
\dots,~$\Delta_N$ be closed disks about $t_1$, \dots,~$t_N$ mutually
disjoint and not containing~$t_0$. For $1\leq i\leq N$, let $\delta_i$~be a
path connecting~$t_0$ to a point~$d_i$ of the boundary~$\partial\Delta_i$
of~$\Delta_i$. Each~$\delta_i$ is required not to meet any of the previous
disks except that the end of~$\delta_i$ touches~$\Delta_i$. Moreover paths
$\delta_1$, \dots,~$\delta_N$ are required not to meet together except at
their origin. For $1\leq i\leq N$, let $\omega_i$~be a loop based at~$d_i$
and running once counter-clockwise around~$\partial\Delta_i$. Finally, for
$1\leq i\leq N$, consider the loops
$\gamma_i=\delta_i*\omega_i*\delta_i^-$, where $*$ designates concatenation
and where $\delta_i^-$~is the path opposite to~$\delta_i$. We denote
by~$\bar\gamma_i$ the homotopy class of~$\gamma_i$ in
$\PP^1-\set{t_1,\dots,t_N}$.
\end{dfn}

We shall obtain the wanted monodromies with the help of some special
isotopies above these loops.

\begin{lem}
\label{t:G}
For $1\leq i\leq N$, there are isotopies
\begin{align*}
 G_i      &\colon\quad\Ll_{t_0}\cap(W-j(V))\times\I\longrightarrow W-j(V),\\
 \hat G_i &\colon\quad\Ll_{t_0}^\sharp\cap(\hat W-\widehat{j(V)})\times\I
  \longrightarrow\hat W-\widehat{j(V)}
\end{align*}
such that
\begin{enumerate}
\renewcommand{\labelenumi}{(\Roman{enumi})}
\item $G_i(x,0)=x$ for any $x\in\Ll_{t_0}\cap(W-j(V))$,
\item $G_i(\,.\,,s)$ is a homeomorphism from
 $\Ll_{t_0}\cap(W-j(V))$ onto
 $\Ll_{\gamma_i(s)}\cap(W-j(V))$ for any $s\in\I$,
\item $G_i(y,s)=y$ for any $y\in \Mm\cap(W-j(V))$ and $s\in\I$,
\item[(\^I)] $\hat G_i(v,0)=v$,
 for any $v\in\Ll_{t_0}^\sharp\cap(\hat W-\widehat{j(V)})$,
\item[(\^I\^I)] $\hat G_i(\,.\,,s)$ is a homeomorphism from
 $\Ll_{t_0}^\sharp\cap(\hat W-\widehat{j(V)})$ onto
 $\Ll_{\gamma_i(s)}^\sharp\cap(\hat W-\widehat{j(V)})$ for any $s\in\I$,
\item[(\^I\^I\^I)] $\hat G_i((y,t_0),s)=(y,\gamma_i(s))$ for any
 $y\in\Mm\cap(W-j(V))$ and $s\in\I$.
\end{enumerate}
Moreover, $G_i$ and~$\hat G_i$ can be asked to fit together by making
commutative the following diagram:
\begin{equation}
\label{GPhi}
 \begin{CD}
 \Ll_{t_0}^\sharp\cap(\hat W-\widehat{j(V)})\times\I
         @>\hat G_i>>   \hat W-\widehat{j(V)}\\
 @VV\Phi\times\id_\I V    @VV\Phi V\\
 \Ll_{t_0}\cap(W-j(V))\times\I
         @>G_i>>        W-j(V).
 \end{CD}
\end{equation}
(recall that the blowing down morphism $\Phi$~induces an
isomorphism between $\Ll_{t_0}^\sharp\cap(\hat W-\widehat{j(V)})$ and
$\Ll_{t_0}\cap(W-j(V))$; following Notation~\ref{t:notn1},
item~\ref{t:notn1.maps}, all maps induced by~$\Phi$ are denoted here by
the same letter).
\end{lem}

\begin{proof}
This follows in a standard manner from Claim~\ref{t:Fibr}, starting
from~$\hat G_i$ and then going down to~$G_i$ thanks to the isomorphism
between $\Ll_t^\sharp\cap(\hat W-\widehat{j(V)})$ and
$\Ll_t\cap(W-j(V))$: cf.~\cite[Lemmas 4.1 and~4.2]{London}. The
statements for points (II) and~(\^I\^I) in~\cite{London} are weaker than
here but the proof given in~\cite{London} is valid for the stronger form.
\end{proof}
It must be noticed that isotopies~$G_i$ and~$\hat G_i$ are not
uniquely determined by loop~$\gamma_i$. But, if diagram~\eqref{GPhi} is
commutative, one of them determines the other.

The ending stage of these isotopies will be the geometric monodromies we
want. Lemma~\ref{t:G} implies the following.

\begin{lem}
\label{t:H}
For $1\leq i\leq N$, we have homeomorphisms
\begin{align*}
 H_i      &\colon\quad\Ll_{t_0}\cap(W-j(V))
  \longrightarrow \Ll_{t_0}\cap(W-j(V)),\\
 \hat H_i &\colon\quad\Ll_{t_0}^\sharp\cap(\hat W-\widehat{j(V)})
  \longrightarrow\Ll_{t_0}^\sharp\cap(\hat W-\widehat{j(V)}),
\end{align*}
defined by setting
\begin{equation*}
 H_i(x)=G_i(x,1)\quad\text{and}\quad\hat H_i(v)=\hat G_i(v,1).
\end{equation*}
These homeomorphisms leave fixed the points of $\Mm\cap(W-j(V))$ and
$\Mm\cap(W-j(V))\times\set{t_0}$ respectively. Moreover, if
diagram~\eqref{GPhi} is commutative, so is the following:
\begin{equation}
\label{HPhi}
 \begin{CD}
 \Ll_{t_0}^\sharp\cap(\hat W-\widehat{j(V)})
         @>\hat H_i>>   \Ll_{t_0}^\sharp\cap(\hat W-\widehat{j(V)})\\
 @VV\vis\Phi V    @VV\vis\Phi V\\
 \Ll_{t_0}\cap(W-j(V))
         @>H_i>>        \Ll_{t_0}\cap(W-j(V)).
 \end{CD}
\end{equation}

\end{lem}

\begin{dfn}
We shall call~$H_i$ a \emph{geometric monodromy of $\Ll_{t_0}\cap(W-j(V))$
relative to $\Mm\cap(W-j(V))$ above~$\gamma_i$}. Similarly for~$\hat H_i$.
\end{dfn}
Notice that, like the isotopies giving rise to them, these geometric
monodromies are not uniquely determined by the choice of loop~$\gamma_i$.
However we have the following invariance property.

\begin{lem}
\label{t:inv}
Given an index~$i$ with $1\leq i\leq N$, another choice of loop~$\gamma_i$
within the same homotopy class~$\bar\gamma_i$ in $\PP^1-\set{t_1,\dots,t_N}$
and another choice of isotopies $G_i$ and~$\hat G_i$ above~$\gamma_i$,
provided they satisfy conditions (I)\nbd (III) and (\^I)\nbd (\^I\^I\^I)
respectively, would lead to geometric monodromies which are isotopic to
$H_i$ and~$\hat H_i$ through isotopies in $\Ll_{t_0}\cap(W-j(V))$ and
$\Ll_{t_0}^\sharp\cap(\hat W-\widehat{j(V)})$ leaving pointwise fixed the
subsets $\Mm\cap(W-j(V))$ and $\Mm\cap(W-j(V))\times \set{t_0}$
respectively. This is true even if the new loop has not the special
form described in Definition~\ref{t:loops}.
\end{lem}

\begin{proof}
Concerning~$\hat G_i$, this is the classical invariance property of
geometric monodromies with an enhancement about fixed points given by the
trivial subbundle of Claim~\ref{t:Fibr}. A similar property holds for~$G_i$
since it can be associated to an isotopy~$\hat G_i$ making
diagram~\eqref{GPhi} commutative.
\end{proof}

Thus, though geometric monodromy~$H_i$ is not
uniquely defined, its isotopy class in $\Ll_{t_0}\cap(W-j(V))$
relative to $\Mm\cap(W-j(V))$ is unique and wholly determined by
homotopy class~$\bar\gamma_i$. This isotopy class can be called \emph{the}
geometric monodromy of $\Ll_{t_0}\cap(W-j(V))$ relative to
$\Mm\cap(W-j(V))$ \emph{associated to~$\bar\gamma_i$}. Similarly
for~$\hat H_i$.

The isomorphisms $H_{i\lstar}$ and~$\bar H_{i\lstar}$
of $H_k(\Ll_{t_0}\cap(W-j(V))$ and
$H_k\bigl(\Ll_{t_0}\cap(W-j(V)),\Mm\cap(W-j(V))\bigr)$
induced by~$H_i$ depend only on this isotopy class. Similarly for the
algebraic monodromies induced by~$\hat H_i$. In particular, to obtain them,
we could use geometric monodromies arising from maps $G_i$ and~$\hat G_i$
satisfying to the looser requirements of Lemma~\ref{t:inv}. For instance
monodromies above the loops~$\Gamma_i$ described before
Definition~\ref{t:loops}. Or even geometric monodromies not satisfying to
the commutativity of diagram~\eqref{HPhi}; the corresponding diagrams at the
homology and relative homology levels would still be commutative.

Nevertheless, for convenience in our forthcoming constructions, we shall
henceforth use special geometric monodromies $H_i$ and~$\hat H_i$ as we have
constructed above, which are associated with a special set of loops as
given in Definition~\ref{t:loops} and which are linked together by the
commutative diagram~\eqref{HPhi}.

\section{The degeneration operator}
\label{s:degener}

In~\cite[section~2]{Annals}, homotopical degeneration operators are
introduced for generic pencils of hyperplane sections of the complement
in~$\CC^{n+1}$ of a hypersurface with isolated singularities (including at
infinity). 

The purpose of this section is to define projective analogs of these for
pencil $(L_t)_{t\in\PP^1}$, acting on the $(n-1)$\nbd th homotopy group of 
each $L_{t_i}-L_{t_i}\cap V$ when $n\geq3$ and on some subgroup of the
fundamental group of each when $n=2$. According to section~\ref{s:homHom},
these homotopy groups are, when $n\geq3$, canonically identified with the
homology groups of $d$\nbd fold covers introduced there and hence it is
enough to define a homological degeneration operator on these homology
groups. This also works when $n=2$ thanks to the morphisms of
Lemma~\ref{t:homHom}, item~\ref{t:homHom.beta2} and still the
isomorphisms of item~\ref{t:homHom.eta}.

\subsection{Homological degeneration operator on the cover}
\label{s:Hdegener}

Suppose~$n\geq2$. For each~$i$, with $1\leq i\leq N$, we define such an 
operator~$D_i$ so that the following diagram is commutative, where
$\Delta_i$~is a disk about~$t_i$ as described in Definition~\ref{t:loops},
and where arrows labeled~$\tau_i$, $\wn_i$ and~$\overline{\Phi_*}$ are to be
defined in the remainder of this section:
\begin{equation}
\label{Deg}
\begin{CD}
  H_n\bigl(P^{-1}(\partial\Delta_i)\cap(\hat W-\widehat{j(V)})\bigr)   
    @>\his{\wn_i}>>  
       H_n\bigl(\Ll_{t_0}^\sharp\cap(\hat W-\widehat{j(V)})\bigr)/
             \im(\hat H_{i\lstar}-\id)\\
  @AA\tau_i A          @VV\vis\overline{\Phi_*} V\\
  H_{n-1}\bigl(\Ll_{t_i}\cap(W-j(V))\bigr)
    @>D_i>>
       H_n\bigl(\Ll_{t_0}\cap(W-j(V))\bigr)/\im(H_{i\lstar}-\id).
\end{CD}
\end{equation}

The arrow labeled~$\overline{\Phi_*}$ is easily defined as follows. The
blowing down morphism~$\Phi$ induces an isomorphism between
$\Ll_{t_0}^\sharp\cap(\hat W-\widehat{j(V)})$ and $\Ll_{t_0}\cap(W-j(V))$
(by~\eqref{Isom}) which gives an isomorphism~$\Phi_*$
between the $n$\nbd th homology groups of these spaces. This in turn
factorizes into an isomorphism~$\overline{\Phi_*}$ as indicated on the
diagram thanks to the commutativity of diagram~\eqref{HPhi}. Recall that, by
the invariance property of Lemma~\ref{t:inv}, isomorphisms~$H_{i\lstar}$
and~$\hat H_{i\lstar}$ depend only on the homotopy class~$\bar\gamma_i$ of
Definition~\ref{t:loops}.

The arrow labeled~$\wn_i$ arises from the Wang sequence of the
fibration of Claim~\ref{t:Fibr} restricted to the part above the
circle~$\partial\Delta_i$. This is detailed below.

The arrow labeled~$\tau_i$ is essentially a tube map in the Poincar\'e
residue sequence for the complement of
$\Ll_{t_i}^\sharp\cap(\hat W-\widehat{j(V)})$ in
$\hat W-\widehat{j(V)}$. Details are also given below.

Operator~$D_i$~depends only on homotopy class~$\bar\gamma_i$. This will be
easy to see after the comparison between the degeneration and variation
operators made in section~\ref{s:link} (see Corollary~\ref{t:Dinv}).

\subsubsection*{Definition of isomorphism~$\wn_i$}

To define~$\wn_i$, we use the Wang sequence of the fibration indicated
above. In this sequence we want to use the fiber
above~$t_0$ of the fibration of Claim~\ref{t:Fibr}, though $t_0$~is outside
of~$\partial\Delta_i$, and the monodromy~$\hat H_i$ above the
loop~$\gamma_i$ of Definition~\ref{t:loops}
instead of the monodromy above~$\partial\Delta_i$.

For this purpose, let us get monodromy~$\hat H_i$ in three steps, following
the decomposition
$\gamma_i=\delta_i*\omega_i*\delta_i^-$ given in Definition~\ref{t:loops}.
Let
$\hat H_i'\colon\Ll_{d_i}^\sharp\cap(\hat W-\widehat{j(V)})
\to\Ll_{d_i}^\sharp\cap(\hat W-\widehat{j(V)})$ be a geometric
monodromy above~$\omega_i$ defined in the same way as~$\hat H_i$,
using Lemmas~\ref{t:H} and~\ref{t:G} but replacing
parameter~$t_0$ by the base point~$d_i$ of~$\omega_i$ and loop~$\gamma_i$
by loop~$\omega_i$ wherever they occur.  Also let
$\hat H_i''\colon\Ll_{t_0}^\sharp\cap(\hat W-\widehat{j(V)})
\to\Ll_{d_i}^\sharp\cap(\hat W-\widehat{j(V)})$ be a
homeomorphism obtained with the same definitions, this time
replacing~$\gamma_i$ by~$\delta_i$. If $\hat H_i'$~is constructed using
isotopy~$\hat G_i'$ and $\hat H_i''$ using~$\hat G_i''$, we can
build~$\hat H_i$ from
$\hat G_i=\hat G_i''*\hat G_i'*\hat G_i^{\p\p-}$, where operations on
isotopies parallel those on paths, each isotopy taking the fiber up in the
place it was left by the former. This is a legal choice for~$\hat G_i$
since it can easily been verified that it satisfies to the conditions of
Lemma~\ref{t:G}. With this setting for~$\hat G_i$, the corresponding
geometric monodromy~$\hat H_i$ is decomposed as
\begin{equation}
\label{hatHH'H''}
 \hat H_i=\hat H_i^{\p\p-1} \circ \hat H_i' \circ \hat H_i''.
\end{equation}

As loop~$\omega_i$ runs once around~$\Delta_i$, the monodromy~$\hat H_i'$
above~$\omega_i$ fits into the Wang sequence of the fibration
above~$\partial\Delta_i$ we consider. We embed this sequence into the
following diagram where it appears as the upper line (we use
Notation~\ref{t:notn2}):
\begin{equation}
\label{Wang}
\begin{CD}
 H_n(\Ll_{d_i}^\shrpp)
   @>\hat H_{i\lstar}'-\id>>  
      H_n(\Ll_{d_i}^\shrpp)    @>\incl_*>>  
          H_n(P^{\p-1}(\partial\Delta_i))    @>>>
                H_{n-1}(\Ll_{d_i}^\shrpp).\\
 @AA\vis\hat H_{i\lstar}''A   
    @AA\vis\hat H_{i\lstar}''A     @.    @.\\
 H_n(\Ll_{t_0}^\shrpp)
  @>\hat H_{i\lstar}-\id>>
     H_n(\Ll_{t_0}^\shrpp)    @.   @.
\end{CD}
\end{equation}
It is commutative by~\eqref{hatHH'H''}. But
$H_{n-1}(\Ll_{d_i}^\shrpp)=0$ because this group is
isomorphic to  $H_{n-1}(\Ll_{d_i}\cap(W-j(V))$ which vanishes when
$n\geq2$ (cf.~\eqref{Isom} and Lemma~\ref{t:vanHom}). Thus the inclusion map
induces an isomorphism
\begin{equation*}
 H_n(\Ll_{d_i}^\shrpp)/\im(\hat H_{i\lstar}'-\id)
   \isarrow H_n(P^{\p-1}(\partial\Delta_i)).
\end{equation*}
Then, by commutativity of diagram~\eqref{Wang}, 
homeomorphism~$\hat H_i''$ followed by inclusion induces also an isomorphism
\begin{equation}
\label{wn^-1}
 H_n(\Ll_{t_0}^\shrpp)/
   \im(\hat H_{i\lstar}-\id)\isarrow
      H_n(P^{\p-1}(\partial\Delta_i)).
\end{equation}
The isomorphism~$\wn_i$ appearing in diagram~\eqref{Deg} is the
inverse of this one.
 
\subsubsection*{Definition of homomorphism~$\tau_i$}

Homomorphism~$\tau_i$ is defined as a tube map in a Poincar\'e residue
sequence (also named Leray or Thom-Gysin sequence) through the following
diagram where we still use Notation~\ref{t:notn2}:
\begin{equation}
\label{tube}
\begin{CD}
 H_{n-1}(\Ll_{t_i}')  @.  @.  \\
 @AA\vis\Phi_* A   @.   @.\\
 H_{n-1}(\Ll_{t_i}^\shrpp)   @>\shis{T_i}>>  
   \mspace{-26mu}H_{n+1}(\hat W',\hat W'-\Ll_{t_i}^\shrpp)    @.\\
 @.     \mspace{-26mu}@AA\vis\incl_* A      @.    \\
 @.  \mspace{-26mu}\makebox[0pt]%
       {$H_{n+1}(P^{\p-1}(\Delta_i),P^{\p-1}
            (\Delta_i)-\Ll_{t_i}^\shrpp)$}    @.\\
 @.  \mspace{-26mu}@AA\vis\incl_* A      @.    \\
 @.   \mspace{-26mu}H_{n+1}(P^{\p-1}(\Delta_i),P^{\p-1}(\partial\Delta_i))
    @>\partial >>    H_n(P^{\p-1}(\partial\Delta_i)).
\end{CD}
\end{equation}

Homomorphism~$\tau_i$ is obtained by overall composition from the upper
left to the lower right end of the diagram (reversing isomorphisms when
necessary). The arrow labeled~$\Phi_*$ is an isomorphism induced by
the blowing down morphism~$\Phi$ (see~\eqref{Isom}).
Following our general convention, arrows labeled~$\incl_*$ are induced by
inclusion. The upper one is an excision isomorphism. The lower one is also
an isomorphism because
$\partial\Delta_i$~is a strong deformation retract of
$\Delta_i-\set{t_i}$ and the spaces $P^{\p-1}(\partial\Delta_i)$ and
$P^{\p-1}(\Delta_i)-\Ll_{t_i}^\shrpp$ are the parts over
$\partial\Delta_i$ and $\Delta_i-\set{t_i}$ of the locally trivial
fibration of Claim~\ref{t:Fibr}, so that the inclusion of the former into
the latter is a homotopy equivalence (see~\cite[proof of Lemme~4.4]{L'Ens}).

The significant arrow is the one labeled~$T_i$ which is a Leray (or
Thom-Gysin) isomorphism. Such an isomorphism arises whenever one removes a
closed submanifold~$P$ from a Hausdorff paracompact complex manifold~$N$.
If $P$~has pure complex codimension~$c$ in~$N$, this is an isomorphism from
$H_{k-2c}(P)$ onto $H_k(N,N-P)$ holding for any~$k$, with the convention
that $H_{k-2c}(P)=0$ for~$k<2c$ (cf.~\cite[Annexe]{L'Ens}). Here we apply it
with $N=\hat W'$, $P=\Ll_{t_i}^\shrpp$, $c=1$ and $k=n+1$. We must verify
that these settings satisfy the above conditions on~$N$, $P$ and~$c$.

First, $\hat W'=\hat W-\widehat{j(V)}$ is the total transform of
$W-j(V)$ when blowing~$\PP^{n+2}$ up (cf.~section~\ref{s:Blow}). But
$W-j(V)$ is a submanifold of~$\PP^{n+2}$ by Claim~\ref{t:Wsing} and the
$n$\nbd plane~$\Mm$ along which $\PP^{n+2}$~is blown up is,
by Claim~\ref{t:Transv}, transverse to $W-j(V)$. It follows that the total
transform~$\hat W'$ of $W-j(V)$ is a submanifold of~$\HPP^{n+2}$
(cf.~\cite[(5.5.1)]{L'Ens}). It is Hausdorff paracompact since
$\HPP^{n+2}$~is a subspace of $\PP^{n+2}\times\PP^1$ which is
metrizable.

Second, we have $\Ll_{t_i}^\shrpp=\Ll_{t_i}^\sharp\cap\hat
W'$ and $\Ll_{t_i}^\sharp$~is the strict transform
of~$\Ll_{t_i}$, a member of the pencil~$\Ll$ with base
locus~$\Mm$ introduced in section~\ref{s:Blow}. As $\Ll_{t_i}$~is,
by Claim~\ref{t:Transv}, also transverse to $W-j(V)$, it follows that
$\Ll_{t_i}^\sharp$~is transverse in~$\HPP^{n+2}$ to the
total transform~$\hat W'$ of $W-j(V)$ (cf.~\cite[(5.5.2)]{L'Ens}). But
$\Ll_{t_i}^\sharp$~is a closed submanifold of~$\HPP^{n+2}$
of pure complex codimension~$1$ as it follows from~\eqref{Strict} and the
fact that $P$~is a submersion. Hence $\Ll_{t_i}^\shrpp$~is a
closed submanifold of~$\hat W'$ of pure complex codimension~$1$. The
conditions of validity of the Leray isomorphism are thus checked.

This completes the definition of homomorphism~$\tau_i$ and hence of the
homological degeneration operator~$D_i$ at the $d$\nbd fold cover level.

\subsection{Homotopical degeneration operator}

For each~$i$, with $1\leq i\leq N$, the isomorphism~$\alpha_{t_0}$
and the homomorphisms~$\beta_i$ of Lemma~\ref{t:homHom} will allow us,
if~$n\geq2$, to define a homotopical degeneration operator~$\Dd_i$ from the
homology operator~$D_i$ constructed at the $d$\nbd fold level in the
preceding subsection.

We first define monodromies on $\pi_n(L_{t_0}-L_{t_0} \cap V,e)$ as the
pull-backs of the monodromies on $H_n\bigl(\Ll_{t_0}\cap(W-j(V))\bigr)$ by
isomorphism~$\alpha_{t_0}$.

\begin{dfn}
\label{t:hihash}
Let $e$~be a base point in $M-M\cap V$ as in Lemma~\ref{t:homHom}.
If~$n\geq2$ and for
$1\leq i\leq N$, monodromy~$h_{i\lhash}$ is defined by the commutativity of
the following diagram:
\begin{equation}
\label{halpha}
\begin{CD}
 H_n\bigl(\Ll_{t_0}\cap(W-j(V))\bigr)
      @>\his{H_{i\lstar}}>>    
                     H_n\bigl(\Ll_{t_0}\cap(W-j(V))\bigr)\\
 @AA\vis\alpha_{t_0}A                       @AA\vis\alpha_{t_0}A\\
 \pi_n(L_{t_0}-L_{t_0}\cap V,e)
      @>\shis{h_{i\lhash}}>>   
                     \pi_n(L_{t_0}-L_{t_0}\cap V,e).
\end{CD}
\end{equation}
\end{dfn}

As $H_{i\lstar}$ depends only on the homotopy class~$\bar\gamma_i$ of
Definition~\ref{t:loops}, so does monodromy~$h_{i\lhash}$.

\begin{rem}
\label{t:truehihash}
In fact $h_{i\lhash}$~is indeed induced on the $n$\nbd th homotopy group by
a geometric monodromy~$h_i$ of $L_{t_0}-L_{t_0}\cap V$ as the notation
suggests. Such a monodromy is obtained in the same way as $H_i$~was defined
from isotopy~$G_i$, by using an isotopy~$g_i$ satisfying to conditions
similar to those given for~$G_i$ in Lemma~\ref{t:G}. This exists
by~\cite[Lemma~4.1]{London}. Then $h_i$~satisfies to an invariance property
similar to that of Lemma~\ref{t:inv} and induces an automorphism of
$\pi_n(L_{t_0}-L_{t_0}\cap V,e)$ depending only on~$\bar\gamma_i$. But
$h_i$ and~$H_i$ can be chosen so that they commute with the covering
projection~$\pi$. To do this, one starts from~$g_i$ and defines~$G_i$ as a
lift of~$g_i$ by~$\pi$ satisfying to the initial condition~(I) of
Lemma~\ref{t:G}. Then one can see that $G_i$~satisfies also automatically to
conditions~(II) and~(III) thanks to the similar conditions satisfied
by~$g_i$. The corresponding geometric monodromies $h_i$ and~$H_i$ commute
with~$\pi$ as desired. This fact together with the functoriality of the
Hurewicz homomorphisms entail that the induced automorphism~$h_{i\lhash}$
of $\pi_n(L_{t_0}-L_{t_0}\cap V,e)$ makes diagram~\eqref{halpha} commutative
(recall the definition of~$\alpha_{t_0}$ in Lemma~\ref{t:homHom}). Hence
this~$h_{i\lhash}$ coincides with the one of Definition~\ref{t:hihash}.
\end{rem}

If $n\geq2$ and for $1\leq i\leq N$, commutative diagram~\eqref{halpha}
allows to define in turn an isomorphism~$\overline{\alpha_{t_0}}$ making
commutative the following diagram:
\begin{equation}
\label{alphaalpha}
\begin{CD}
 H_n\bigl(\Ll_{t_0}\cap(W-j(V))\bigr)    @>\can>>  
         H_n\bigl(\Ll_{t_0}\cap(W-j(V))\bigr)/\im(H_{i\lstar}-\id)\\
 @AA\vis\alpha_{t_0}A           @AA\vis\overline{\alpha_{t_0}}A\\
 \pi_n(L_{t_0}-L_{t_0} \cap V,e)         @>\can>>
         \pi_n(L_{t_0}-L_{t_0}\cap V,e)/\im(h_{i\lhash}-\id).
\end{CD}
\end{equation}

Then, if $n\geq3$ and for $1\leq i\leq N$, isomorphism~$\beta_i$ of
Lemma~\ref{t:homHom}, item~\ref{t:homHom.beta3} together with this
isomorphism~$\overline{\alpha_{t_0}}$ lead from homological operator~$D_i$
to the homotopical degeneration operator~$\Dd_i$ we are looking for. This
is done by asking the following diagram to be commutative:
\begin{equation}
\label{betaalpha}
\begin{CD}
 H_{n-1}\bigl(\Ll_{t_i}\cap(W-j(V))\bigr)    @>D_i>>  
         H_n\bigl(\Ll_{t_0}\cap(W-j(V))\bigr)/\im(H_{i\lstar}-\id)\\
 @AA\vis\beta_i A           @AA\vis\overline{\alpha_{t_0}}A\\
 \pi_{n-1}(L_{t_i}-L_{t_i}\cap V,e)          @>\Dd_i>>
         \pi_n(L_{t_0}-L_{t_0}\cap V,e)/\im(h_{i\lhash}-\id).
\end{CD}
\end{equation}
When $n=2$, isomorphism~$\overline{\alpha_{t_0}}$ and this time
homomorphism~$\beta_i$ of Lemma~\ref{t:homHom}, item~\ref{t:homHom.beta2}
lead to an operator~$\Dd_i$ defined on the subgroup~$G_i$ introduced there,
by asking the following diagram to be commutative:
\begin{equation}
\label{betaalpha2}
\begin{CD}
 H_1\bigl(\Ll_{t_i}\cap(W-j(V))\bigr)    @>D_i>>  
         H_2\bigl(\Ll_{t_0}\cap(W-j(V))\bigr)/\im(H_{i\lstar}-\id)\\
 @AA\beta_i A           @AA\vis\overline{\alpha_{t_0}}A\\
 G_i          @>\Dd_i>>
         \pi_2(L_{t_0}-L_{t_0}\cap V,e)/\im(h_{i\lhash}-\id).
\end{CD}
\end{equation}
Like~$D_i$, operator~$\Dd_i$ depends only on the homotopy
class~$\bar\gamma_i$ of Definition~\ref{t:loops}.

\section{The variation operator}
\label{s:var}

In \cite[section~4]{London} homological variation operators are defined for
generic pencils of hyperplane sections of a quasi-projective variety. They
are analogous to the classical variation operator associated with the
Milnor fibration of an isolated singularity
(see~\cite[chapter~2]{Singularity}). In our situation they give homological
variation operators for pencil
$(L_t)_{t\in\PP^1}$, defined on the $n$\nbd th relative homology group of
$L_{t_0}-L_{t_0}\cap V$ modulo
$M-M\cap V$ and associated with each special member~$L_{t_i}$ of the
pencil, more precisely with the homotopy class~$\bar\gamma_i$ in
$\PP^1-\set{t_1,\dots,t_N}$ of a loop~$\gamma_i$ surrounding~$t_i$ in the
parameter space as in Definition~\ref{t:loops}.

In this section we want to define, when~$n\geq2$, homotopical analogs of
these,
\begin{equation*}
 \Vv_i\colon \pi_n(L_{t_0}-L_{t_0}\cap V,M-M\cap V,e)
 \longrightarrow \pi_n(L_{t_0}-L_{t_0}\cap V,e)
\end{equation*}
associated with~$\bar\gamma_i$ for $1\leq i\leq N$, where $e$~is a base
point in $M-M\cap V$.

As for the degeneration operators, we shall go to the $d$\nbd fold cover
level and use homological variation operators defined there,
\begin{equation*}
 V_i\colon H_n\bigl(\Ll_{t_0}\cap(W-j(V)), \Mm\cap(W-j(V))\bigr)
 \longrightarrow H_n\bigl(\Ll_{t_0}\cap(W-j(V))\bigr)
\end{equation*}
associated, for $1\leq i\leq N$, with the homotopy classes~$\bar\gamma_i$
of Definition~\ref{t:loops}.

We recall the definition and properties of operator~$V_i$ as given
in~\cite[section~4]{London}, which in fact hold with~$n\geq1$.

For any relative $n$\nbd cycle~$\Xi$ on $\Ll_{t_0}\cap(W-j(V))$ modulo
$\Mm\cap(W-j(V))$, one defines
\begin{equation}
\label{Vardef}
 V_i(\class{\Xi}_{(\Ll_{t_0}\cap(W-j(V)),\Mm\cap(W-j(V)))})
 =\class{H_{i\lbullet}(\Xi)-\Xi}_{\Ll_{t_0}\cap(W-j(V))}
\end{equation}
using Notation~\ref{t:notn1}, items \ref{t:notn1.topalg}
and~\ref{t:notn1.Hom}. Due to the fact that $H_i$~leaves the points of
$\Mm\cap(W-j(V))$ fixed (Lemma~\ref{t:H}), the chain
$H_{i\lbullet}(\Xi)-\Xi$ is actually an absolute cycle and the
correspondence $\Xi\mapsto H_{i\lbullet}(\Xi)-\Xi$ induces a
homomorphism~$V_i$ at the homology level (\cite[Lemmas 4.6
and~4.8]{London}). Thanks to the invariance property expressed by
Lemma~\ref{t:inv}, this homomorphism depends only on homotopy
class~$\bar\gamma_i$ (\cite[Lemma~4.8]{London}).

Now, if $n\geq2$ and for $1\leq i\leq N$, isomorphism~$\alpha_{t_0}$ of
Lemma~\ref{t:homHom} and homomorphism~$\bar\alpha_{t_0}$ of
Lemma~\ref{t:relhomHom} lead from~$V_i$ to the wanted homotopical
variation operator~$\Vv_i$ by asking the following diagram to be
commutative:
\begin{equation}
\label{Valpha}
\begin{CD}
 H_n\bigl(\Ll_{t_0}\cap(W -j(V)),\Mm\cap(W-j(V))\bigr)    @>V_i>>
                              H_n\bigl(\Ll_{t_0}\cap(W-j(V))\bigr) \\ 
 @AA\bar\alpha_{t_0}A   @AA\vis\alpha_{t_0}A \\ 
 \pi_n(L_{t_0}-L_{t_0}\cap V,M-M\cap V,e)    @>\Vv_i>>
                              \pi_n (L_{t_0}-L_{t_0}\cap V,e). \\ 
\end{CD}
\end{equation}
As $V_i$ depends only on the homotopy class~$\bar\gamma_i$ of
Definition~\ref{t:loops}, so does operator~$\Vv_i$.

\begin{rem}
The homological variation operators
\begin{equation*}
 v_i\colon H_n(L_{t_0}-L_{t_0}\cap V,M-M\cap V)
 \longrightarrow H_n(L_{t_0}-L_{t_0}\cap V)
\end{equation*}
we talked about at the beginning of the section are given by a formula
similar to~\ref{Vardef} using the monodromies~$h_i$ considered in
Remark~\ref{t:truehihash}. The homotopical variation operators~$\Vv_i$ we
have defined here are linked to those by Hurewicz homomorphisms as is shown
in the following diagram:
\begin{equation*}
\begin{CD}
 \pi_n(L_{t_0}-L_{t_0}\cap V,M-M\cap V,e)    @>\Vv_i>>
                              \pi_n (L_{t_0}-L_{t_0}\cap V,e) \\ 
  @VV\bar\chi V   @VV\chi V \\
 H_n(L_{t_0}-L_{t_0}\cap V,M-M\cap V)   @>v_i>>
                               H_n (L_{t_0}-L_{t_0}\cap V). \\
\end{CD}
\end{equation*}
The commutativity of this diagram is a consequence of the commutativity of
the diagram~\eqref{Valpha} defining~$\Vv_i$, of the definitions of
homomorphisms $\alpha_{t_0}$ and~$\bar\alpha_{t_0}$ occurring there (see
Lemmas \ref{t:homHom} and~\ref{t:relhomHom}), of the functoriality of
Hurewicz homomorphisms and of the commutation of monodromies $h_i$
and~$H_i$ with the covering projection~$\pi$ as stated in
Remark~\ref{t:truehihash}.
\end{rem}

We end this section by noticing that the homotopical variation
operator~$\Vv_i$ when restricted to absolute cycles acts like the variation
$h_{i\lhash}-\id$ of the homotopical monodromy associated
with~$\bar\gamma_i$. This is specified by the following lemma.

\begin{lem}
\label{t:varmonodr}
If $n\geq2$ then, for $1\leq i\leq N$, the following diagram is
commutative:
\begin{equation}
\label{h-id.var}
\begin{CD}
 \pi_n(L_{t_0}-L_{t_0}\cap V,e)   @>\incl_\#>>
     \pi_n(L_{t_0}-L_{t_0}\cap V,M-M\cap V,e)\\
 @.            @VV{\Vv_i}V\\
 \darrow{65}{17}{-17}{97}{0}{6}{h_{i\lhash}-\id}
    @.            \pi_n(L_{t_0}-L_{t_0}\cap V,e).
\end{CD}
\end{equation}
\end{lem}
\begin{proof}
This will follow from the commutativity of the corresponding homology
diagram at the $d$\nbd fold covering level:
\begin{equation}
\label{H-id.Var}
\begin{CD}
 H_n\bigl(\Ll_{t_0}\cap(W-j(V))\bigr)   @>\incl_*>>
     H_n\bigl(\Ll_{t_0}\cap(W-j(V)),\Mm\cap(W-j(V))\bigr)\\
 @.            @VV{V_i}V\\
 \darrow{78}{17}{-15}{115}{0}{6}{H_{i\lstar}-\id}
    @.            H_n\bigl(\Ll_{t_0}\cap(W-j(V))\bigr).
\end{CD}
\end{equation}
Indeed, diagrams \eqref{h-id.var} and
\eqref{H-id.Var} are linked together by homomorphisms $\alpha_{t_0}$
and~$\bar\alpha_{t_0}$ and one can use the commutativity of diagrams
\eqref{alphagamma}, \eqref{Valpha} and~\eqref{halpha} and the injectivity
of~$\alpha_{t_0}$. The commutativity of diagram~\eqref{H-id.Var} can be
checked in turn by a straightforward computation at the chain level using
formula~\eqref{Vardef}.
\end{proof}

\section{The link between degeneration and variation operators}
\label{s:link}

In this section we make the link between the degeneration
operators~$\Dd_i$ defined in section~\ref{s:degener} and the variation
operators~$\Vv_i$ defined in section~\ref{s:var}. As a side result, we shall
obtain the invariance property for degeneration operators~$D_i$ and
hence~$\Dd_i$, stated in section~\ref{s:degener}.

The main result is the following.

\begin{prop}
\label{t:vardeg}
Using Notation~\ref{t:notn2}, we have, for $1\leq i\leq N$,
the following commutative diagram if~$n\geq3$:
\begin{equation*}
\begin{CD}
\pi_{n-1}(M',e)    @>\incl_\#>>   \pi_{n-1}(L_{t_i}',e)   
      @>\Dd_i>>   \pi_n(L_{t_0}',e)/\im(h_{i\lhash}-\id)\\
@AA\partial A      @.          @AA\can A\\
\pi_n(L_{t_0}',M',e)   @.
    \harrow{20}{0}{135}{\Vv_i}
         @.     \pi_n(L_{t_0}',e).
\end{CD}
\end{equation*}
If $n=2$, then:
\begin{equation*}
\begin{CD}
H    @>\incl_\#>>   G_i   
      @>\Dd_i>>   \pi_2(L_{t_0}',e)/\im(h_{i\lhash}-\id)\\
@AA\partial A      @.          @AA\can A\\
\pi_2(L_{t_0}',M',e)   @.
    \harrow{20}{0}{100}{\Vv_i}
         @.     \pi_2(L_{t_0}',e),
\end{CD}
\end{equation*}
where groups $H$ and~$G_i$~are defined in Lemma~\ref{t:homHom},
item~\ref{t:homHom.beta2} and homomorphisms $\incl_\#$ and~$\partial$ are
well defined by the same reference and Lemma~\ref{t:relhomHom}.
\end{prop}

Before proving this result, we state a corollary relating the images of
the considered operators.

\begin{coro}
\label{t:imvardeg}
If $n\geq2$, then, for $1\leq i\leq N$,
\begin{equation*}
\im\Dd_i=\im\Vv_i\,/\im(h_{i\lhash}-\id).
\end{equation*}
\end{coro}

This makes sense since $\Dd_i$~takes its values in
$\pi_n(L_{t_0}-L_{t_0}\cap V,e)/\im(h_{i\lhash}-\id)$ while
$\Vv_i$~takes its values in $\pi_n(L_{t_0}-L_{t_0}\cap V,e)$ with an
image containing $\im(h_{i\lhash}-\id)$ by
Lemma~\ref{t:varmonodr}.

\begin{proof}[Proof of Corollary~\ref{t:imvardeg}]
The inclusion
$\im\Dd_i\supset\im\Vv_i\,/\im(h_{i\lhash}-\id)$ is clear from
the diagrams of Proposition~\ref{t:vardeg}. The reverse inclusion will also
be clear from it, once proved the following lemma.
\end{proof}

\begin{lem}
\label{t:surj}
The homomorphisms $\partial$ and $\incl_\#$ in the diagrams of
Proposition~\ref{t:vardeg} are surjective.
\end{lem}

\begin{proof}
The case $n=2$ forces us to go into the $d$\nbd fold covers. It is then
more economical to treat the general case thus. Covering projection~$\pi$
induces an isomorphism from $\pi_n(\Ll_{t_0}',\Mm',\varepsilon)$ onto
$\pi_n(L_{t_0}',M',e)$ (see Lemma~\ref{t:relhomHom}) and isomorphisms from
$\pi_{n-1}(\Mm',\varepsilon)$ and $\pi_{n-1}(\Ll_{t_i}',\varepsilon)$ onto
$\pi_{n-1}(M',e)$ when $n\geq3$ (resp.~$H$ when $n=2$) and
$\pi_{n-1}(L_{t_i}',e)$ (resp.~$G_i$ when $n=2$) (see Lemma~\ref{t:homHom}).
It will then be enough to prove that homomorphisms
$\partial\colon\pi_n(\Ll_{t_0}',\Mm',\varepsilon)
\to\pi_{n-1}(\Mm',\varepsilon)$ and
$\incl_\#\colon\pi_{n-1}(\Mm',\varepsilon)
\to\pi_{n-1}(\Ll_{t_i}',\varepsilon)$ are surjective. This is the case for
homomorphism~$\partial$ due to the homotopy exact sequence of the pair
$(\Ll_{t_0}',\Mm')$ and to the fact that
$\pi_{n-1}(\Ll_{t_0}',\varepsilon)$ is trivial if $n\geq3$ as noticed in
the proof of Lemma~\ref{t:homHom} and also when $n=2$ by
Lemma~\ref{t:univcov}. As to homomorphism~$\incl_\#$, it is surjective due
to the Lefschetz hyperplane section theorem for non-singular
quasi-projective varieties (cf.~\cite[1.1.3]{H-L} or~\cite[II.5.1]{G-M2})
when applied to hyperplane~$\Mm\subset\Ll_{t_i}$ cutting the
quasi-projective variety~$\Ll_{t_i}'=\Ll_{t_i}\cap(W-j(V))$ which is
non-singular and of pure dimension~$n$ by Claim~\ref{t:Wsing} and
Claim~\ref{t:Transv}. For the validity of the quoted theorem,
hyperplane~$\Mm$ must fulfill some condition of genericity. By~\cite[Lemma
of the Appendix]{festschrift}, or~\cite[the remark ending the proof
of~II.5.1]{G-M2}, it is enough that
$\Mm$~be transverse to all the strata of a Whitney stratification of
$\Ll_{t_i}\cap W$ having $\Ll_{t_i}\cap j(V)$ as a union of strata. But,
thanks to the transversality of~$\Mm$ to the strata of the
stratification~$\Sigma$ of~$W$ defined in Claim~\ref{t:Strat}, the trace
of~$\Sigma$ on~$\Ll_{t_i}$ can be refined into such a stratification of
$\Ll_{t_i}\cap W$ (cf.~\cite[lemme~11.3]{L'Ens}).
\end{proof}

The sequel of this section will be devoted to the proof of
Proposition~\ref{t:vardeg} and to the result on the invariance of
operator~$D_i$ mentioned above.

\subsubsection*{Proof of Proposition~\ref{t:vardeg}}

We shall go up to the homology of $d$\nbd fold coverings and consider
for $n\geq2$ the following homology diagram which corresponds at this level
to the diagrams of Proposition~\ref{t:vardeg}. It also uses
Notation~\ref{t:notn2}.
\begin{equation}
\label{VarDeg}
\begin{CD}
 H_{n-1}(\Mm')    @>\incl_*>>  H_{n-1}(\Ll_{t_i}')   
          @>D_i>>   H_n(\Ll_{t_0}')/\im(H_{i\lstar}-\id)\\
 @AA\partial A      @.            @AA\can A\\
 H_n(\Ll_{t_0}',\Mm')   @.
    \harrow{19}{0}{130}{V_i}
          @.     H_n(\Ll_{t_0}').
\end{CD}
\end{equation}
It will be enough to show that this diagram is commutative
since it is linked to the first or second
diagram of Proposition~\ref{t:vardeg}, according to whether $n\geq3$ or
$n=2$, by the commutative diagrams \eqref{alphagamma} or \eqref{alphagamma2}
(right-hand parts), \eqref{gammabeta} or \eqref{gammabeta2},
\eqref{betaalpha} or \eqref{betaalpha2}, \eqref{Valpha}
and~\eqref{alphaalpha}, and since the homomorphism~$\overline{\alpha_{t_0}}$
(defined by~\eqref{alphaalpha}) linking the upper right corners of these
two diagrams is injective.

Before proving the commutativity of diagram~\ref{VarDeg}, we first
treat the invariance property for
operators~$D_i$ mentioned in section~\ref{s:Hdegener}.

\begin{coro}
\label{t:Dinv}
For~$n\geq2$ and $1\leq i\leq N$, operator~$D_i$ depends only
on the homotopy class~$\bar\gamma_i$
of Definition~\ref{t:loops}.
\end{coro}

\begin{proof}
As we saw in section~\ref{s:var}, the same invariance property holds for
operator~$V_i$. Then, the commutativity of diagram~\eqref{VarDeg} implies
the assertion for~$D_i$ because, as we shall see, the arrows
labeled~$\partial$ and~$\incl_*$ are surjective. The surjectivity of
homomorphism~$\partial$ results from the homology exact sequence of the pair
$(\Ll_{t_0}',\Mm')$ and the fact that
$H_{n-1}(\Ll_{t_0}')=0$ if~$n\geq2$ by Lemma~\ref{t:vanHom}. The
homomorphism induced by inclusion~$\incl_*$ is surjective due to the
Lefschetz hyperplane section theorem for non-singular quasi-projective
varieties with the same justification as in the proof of Lemma~\ref{t:surj}
but this time applied to homology.
\end{proof}

The commutativity of diagram~\ref{VarDeg} is a consequence of the
following. On one hand, the bundle of Claim~\ref{t:Fibr} has a
trivial subbundle preserved by the isotopy built in section~\ref{s:monodr}
which has a trivial form on it, allowing thus the very definition of the
homological variation operators. On the other hand, this subbundle extends
to a product $\Mm'\times\PP^1$ which is transverse to
each~$\Ll_{t_i}^\shrpp$, so that the tube maps entering in the definition
of the homological degeneration operators have also a trivial form when
restricted to $\Mm'\times\PP^1$. These two facts lead to the link between
operators~$V_i$ and~$D_i$ expressed by the commutativity of
diagram~\ref{VarDeg}. In fact similar considerations come up in the proof
of Proposition~4.13 of~\cite{London} and our assertion will be obtained
along the same lines.

More precisely, we imbed diagram~\eqref{VarDeg} in a larger one, putting the
following diagram on its top (we still use Notation~\ref{t:notn2}):
\begin{equation}
\label{DegWang}
\begin{CD}
 H_n(\Mm'\times\partial\Delta_i)    @>\incl_*>>  
   H_n(P^{\p-1}(\partial\Delta_i))
      @>\his{\wn_i}>>  
         H_n(\Ll_{t_0}^\shrpp)/\im(\hat H_{i\lstar}-\id)\\
 @AA\kappa_i A      @AA\tau_i A          @VV\vis\overline{\Phi_*} V\\
 H_{n-1}(\Mm')    @>\incl_*>>
    H_{n-1}(\Ll_{t_i}')    @>D_i>>
          H_n(\Ll_{t_0}')/\im(H_{i\lstar}-\id).
\end{CD}
\end{equation}
The right-hand square is just diagram~\eqref{Deg}, which defines~$D_i$, and
homomorphism $\kappa_i$~is given by the formula
\begin{equation}
\label{kappa}
 \kappa_i(z)=(-1)^{n-1}z\times\class{\omega_i}_{\partial\Delta_i}\qquad
 \text{for $z\in H_{n-1}(\Mm')$,}
\end{equation}
using the cross-product by the fundamental
class~$\class{\omega_i}_{\partial\Delta_i}$ of~$\partial\Delta_i$, the
loop~$\omega_i$ of Definition~\ref{t:loops} being this time considered as a
$1$\nbd cycle.

Now, to prove the commutativity of diagram~\ref{VarDeg}, it will be enough
to show that diagram~\eqref{DegWang} is commutative as well as the
following diagram which is the outer square of the big diagram obtained by
putting diagrams \eqref{VarDeg} and~\eqref{DegWang} on top of each other:
\begin{equation}
\label{Outer}
\begin{CD}
 H_n(\Mm'\times\partial\Delta_i)    @>\incl_*>>  
    H_n(P^{\p-1}(\partial\Delta_i))
       @>\his{\wn_i}>>  
         H_n(\Ll_{t_0}^\shrpp)/\im(\hat H_{i\lstar}-\id)\\
 @AA\kappa_i A      @.          @VV\vis\overline{\Phi_*} V\\
 H_{n-1}(\Mm')    @.  @.
          H_n(\Ll_{t_0}')/\im(H_{i\lstar}-\id)\\
 @AA\partial A     @.            @AA\can A\\
 H_n(\Ll_{t_0}',\Mm')   @.
    \harrow{14}{0}{165}{V_i}
          @.     H_n(\Ll_{t_0}').
\end{CD}
\end{equation}

\paragraph{Proof of the commutativity of diagram~\eqref{DegWang}}

The right-hand square of this diagram is commutative since it is the
commutative diagram~\eqref{Deg} defining~$D_i$. To see the commutativity of
the left-hand square, let us consider again diagram~\eqref{tube} through
which homomorphism~$\tau_i$ was defined. There is an analogous diagram
obtained by restricting to~$\Mm$ or~$\hat\Mm$ in order to take advantage
of the product structure $\hat\Mm=\Mm\times\PP^1$. Here is
this diagram:
\begin{equation}
\label{tube'}
\begin{CD}
 H_{n-1}(\Mm')  @.  @.  \\
 @AA\vis\Phi_* A   @.   @.\\
 H_{n-1}(\Mm'\times\set{t_i})   @>\shis{T_i'}>>  
    \mspace{66mu}\makebox[0pt]%
             {$H_{n+1}\bigl(\Mm'\times\PP^1,
                \Mm'\times(\PP^1-\set{t_i})\bigr)$}      @.\\
 @.     \mspace{18mu}@AA\vis\incl_* A      @.    \\
 @. \mspace{66mu}\makebox[0pt]%
             {$H_{n+1}\bigl(\Mm'\times\Delta_i,
                \Mm'\times(\Delta_i-\set{t_i})\bigr)$}      @.\\
 @.      \mspace{18mu}@AA\vis\incl_* A      @.    \\
 @.   H_{n+1}(\Mm'\times\Delta_i,\Mm'\times\partial\Delta_i)
    @>\partial >>    H_n(\Mm'\times\partial\Delta_i).
\end{CD}
\end{equation}
As in dagram~\eqref{tube}, the top isomorphism is induced by the blowing
down morphism~$\Phi$ (see~\eqref{Prod}), the upper arrow labeled~$\incl_*$
is an excision isomorphism and the lower one is an isomorphism since
$\partial\Delta_i$~is a deformation retract of~$\Delta_i-\set{t_i}$. The
arrow labeled~$T_i'$ is again a Leray isomorphism. The conditions of
validity would be easy to check directly but they will also follow from a
naturality property we shall consider in a moment.

Each space occurring in diagram~\eqref{tube'} is contained in the
corresponding space of diagram~\eqref{tube}. This is clear for
$\Mm'\subset\Ll_{t_i}'$ since
$\Mm\subset\Ll_{t_i}$ and it can be seen, using~\eqref{Prod}
and~\eqref{Strict}, that all other spaces of~\eqref{tube'} are
the intersections of the corresponding spaces of~\eqref{tube} with
$\hat\Mm=\Mm\times\PP^1$. Thus diagram~\eqref{tube'} is linked
to diagram~\eqref{tube} by homomorphisms induced by inclusions. All
resulting squares are commutative. This simply follows from the
commutativity of the corresponding diagrams of maps or from the
functoriality of the boundary homomorphism, except for the commutativity of
the following diagram which deserves to be commented on:
\begin{equation}
\label{TT'}
\begin{CD}
 H_{n-1}(\Ll_{t_i}^\shrpp)      @>\shis{T_i}>>  
   \mspace{-26mu}H_{n+1}(\hat W',\hat W'-\Ll_{t_i}^\shrpp)\\
 @AA\incl_*A   @AA\incl_*A\\
 H_{n-1}(\Mm'\times\set{t_i})   @>\shis{T_i'}>>
   H_{n+1}\bigl(\Mm'\times\PP^1,
                \Mm'\times(\PP^1-\set{t_i})\bigr).
\end{CD}
\end{equation}

The commmutativity of this diagram results from the following naturality
property for the Leray isomorphism. With the same notation and hypotheses as
in the exposition we gave of it two paragraphs after diagram~\eqref{tube},
suppose that $N'$~is a closed complex submanifold of~$N$ transverse to~$P$
and let
$P'=N'\cap P$. Then the validity conditions are also satisfied for a Leray
isomorphism from $H_{k-2c}(P')$ onto $H_k(N', N'-P')$ and the diagram formed
by the two Leray isomorphisms and the homomorphisms induced by inclusions is
commutative (cf.~\cite[Annexe]{L'Ens}). Applying these facts with
$N=\hat W'$, $P=\Ll_{t_i}^\shrpp$, $c=1$, $k=n+1$ as before and
$N'=\Mm'\times\PP^1$, we find diagram~\eqref{TT'} since
$\Mm'\times\set{t_i}= (\Mm'\times\PP^1)\cap\Ll_{t_i}^\shrpp$ (still
by~\eqref{Prod} and~\eqref{Strict}. But we must verify that this setting
for~$N'$ satisfies the conditions above.

Let us come back to the third paragraph after diagram~\eqref{tube} where we
checked the conditions of validity of the Leray isomorphism~$T_i$. The
properties which gave us that $\hat W'$~is a submanifold of~$\HPP^{n+2}$
give also, by the same reference, that $\hat W'$~is transverse to~$\hat\Mm$
in~$\HPP^{n+2}$. Hence $\hat\Mm\cap\hat W'=\Mm'\times\PP^1$ is a submanifold
of~$\hat W'$. It is closed since $\hat\Mm$~is closed in~$\HPP^{n+2}$. Next,
the properties which gave us that $\Ll_{t_i}^\sharp$~is transverse to~$\hat
W'$ in~$\HPP^{n+2}$ give in fact, by the same reference, that
$\Ll_{t_i}^\sharp$~is transverse to $\hat\Mm\cap\hat W'$ in~$\HPP^{n+2}$. As
$\hat W'$~contains $\hat\Mm\cap\hat W'$, it follows that
$\hat\Mm\cap\hat W'=\Mm'\times\PP^1$ is transverse to
$\Ll_{t_i}^\sharp\cap\hat W'=\Ll_{t_i}^\shrpp$ in~$\hat W'$
(cf.~\cite[proof of Lemme~9.2~(iii)]{L'Ens}). The conditions for the
natural behavior of the Leray isomorphisms~$T_i$ and~$T_i'$ are thus
checked and the commutativity of diagram~\eqref{TT'} is proved.

Thus diagrams~\eqref{tube'} and~\eqref{tube} are linked in a commutative
diagram by homomorphisms induced by inclusions. Let $\tau_i'$~be the
homomorphism obtained by overall composition from the upper left to the
lower right end of diagram~\eqref{tube'}. As homomorphism~$\tau_i$~is
obtained in the same manner in diagram~\eqref{tube}, we get the
commutativity of the left square of diagram~\eqref{DegWang} but with
$\kappa_i$~replaced by~$\tau_i'$. The commutativity of the original diagram
will then follow from the next lemma.
\begin{lem}
The homomorphism~$\tau_i'$ defined above is equal to the
homomorphism~$\kappa_i$ defined in~\eqref{kappa}.
\end{lem}

\begin{proof}
This is due to the product structure of the spaces in
diagram~\eqref{tube'}, especially to the behavior of the Leray isomorphism
in such a case. With the same notation as in the presentation of this
isomorphism two paragraphs after diagram~\eqref{tube}, suppose that
$N=Q\times R$ where $Q$~and $R$~are complex Hausdorff paracompact manifolds
with $R$~of pure complex dimension~$c$ and suppose that $P=Q\times\set{r}$
with~$r\in R$. Then the conditions of validity hold for a Leray
isomorphism~$T$ from $H_{k-2c}(Q\times\set{r})$ onto $H_k(Q\times R, Q\times
R-Q\times\set{r})$ and this isomorphism takes the following special form. If
$z^\sharp\in H_{k-2c}(Q\times\set{r})$ corresponds to $z\in H_{k-2c}(Q)$ by
the canonical identification of~$Q$ to $Q\times\set{r}$, then
$T(z^\sharp)=z\times w$ where $w\in H_{2c}(R,R-\set{r})$ is the fundamental
class defining the canonical orientation of~$R$ about~$r$
(cf.~\cite[Annexe]{L'Ens}). Here we are in this special case for~$T_i'$,
with $Q=\Mm'$, $R=\PP^1$, $r=t_i$, $c=1$ and $k=n+1$. Remember indeed that
$\Mm$~is transverse to $W-j(V)$ in~$\PP^{n+2}$ so that
$\Mm'=\Mm\cap(W-j(V))~$ is a submanifold of~$\PP^{n+2}$. Thus the Leray
isomorphism~$T_i'$ has the explicit expression
\begin{equation}
\label{T'}
 T_i'(z^\sharp)=z\times u_i  \qquad\text{for $z\in H_{n-1}(\Mm')$}
\end{equation}
where $z^\sharp$~corresponds to~$z$ by the canonical
identification of~$\Mm'$ to $\Mm'\times\set{t_i}$ and where
$u_i\in H_2(\PP^1,\PP^1-\set{t_i})$ is the fundamental class defining the
canonical orientation of~$\PP^1$ about~$t_i$.

Now $\Phi_*(z^\sharp)=z$ by the remark following~\eqref{Prod}. Besides, we
can take a representative~$\Omega_i$ of~$u_i$ which is a relative $2$\nbd
cycle of~$\Delta_i$ modulo~$\partial\Delta_i$. Then, by functoriality of the
cross-product, the composition of the two $\incl_*$~isomorphisms of
diagram~\eqref{tube'} gives an isomorphism, we still denote by~$\incl_*$,
such that
\begin{equation*}
 \incl_*(z\times\class{\Omega_i}_{(\Delta_i,\partial\Delta_i)})=z\times u_i.
\end{equation*}
Combining these facts with~\eqref{T'}, we find that, for~$z\in
H_{n-1}(\Mm')$,
\begin{equation*}
 \tau_i'(z)=\partial(z\times\class{\Omega_i}_{(\Delta_i,\partial\Delta_i)})
 =(-1)^{n-1}z\times\partial\class{\Omega_i}_{(\Delta_i,\partial\Delta_i)}.
\end{equation*}
But, by the special choice of~$\omega_i$ in Definition~\ref{t:loops},
\begin{equation*}
 \partial\class{\Omega_i}_{(\Delta_i,\partial\Delta_i)}
 =\class{\omega_i}_{\partial\Delta_i}
\end{equation*}
and the equality $\tau_i'=\kappa_i$ follows.
\end{proof}

This concludes the proof of the commutativity of the left part and hence of
the whole of diagram~\eqref{DegWang}.

\paragraph{Proof of the commutativity of diagram~\eqref{Outer}}

It will be convenient to concentrate our work to the bundle
$P^{\p-1}(\partial\Delta_i)$ which already was used to define
isomorphism~$\wn_i$ in section~\ref{s:Hdegener}. We hence shall express
variation~$V_i$ by means of a variation operator~$V_i'$ above the
loop~$\omega_i$ of Definition~\ref{t:loops} which runs once
counter-clockwise around~$\partial\Delta_i$. Recal that $d_i$~is the
base point of~$\omega_i$. We still use Notation~\ref{t:notn2}.

\begin{dfn} We define a homological variation operator
\begin{equation*}
	V_i'\colon\quad H_n(\Ll_{d_i}',\Mm')\longrightarrow
	H_n(\Ll_{d_i}')
\end{equation*}
in the same way as~$V_i$ was defined in formula~\eqref{Vardef} but
replacing $\Ll_{t_0}'$ by $\Ll_{d_i}'$ and monodromy~$H_i$ by a
monodromy~$H_i'$ above~$\omega_i$.
\end{dfn}

Just as $V_i$~depends only on the homotopy class of~$\gamma_i$ in
$\PP^1-\set{t_1,\dots,t_N}$, operator~$V_i'$ depends only on the homotopy
class of~$\omega_i$. Therefore operator~$V_i'$ is specified by the
requirement that $\omega_i$~runs once counter-clockwise
around~$\partial\Delta_i$. The monodromy $H_i'\colon\Ll_{d_i}'\to\Ll_{d_i}'$
must of course be defined in the same way as monodromy~$H_i$, using Lemmas
\ref{t:H} and~\ref{t:G} but replacing parameter~$t_0$ by the base
point~$d_i$ of~$\omega_i$ and loop~$\gamma_i$ by loop~$\omega_i$ wherever
they occur (just as we did for the monodromy~$\hat H_i'$ at the blow up
level in the definition of isomorphism~$\wn_i$ in section~\ref{s:Hdegener}).
It will be moreover convenient to have $H_i'$ and~$\hat H_i'$ linked
together by the analog of diagram~\eqref{HPhi}. This is obtained by asking
the commutativity of the analog of diagram~\eqref{GPhi} in Lemma~\ref{t:G}
when building the isotopies leading to $H_i'$ and~$\hat H_i'$. Now, to make
the link with~$V_i$, we choose the monodromy~$H_i$ defining~$V_i$ by
following the same process as we did for~$\hat H_i'$ in the definition of
isomorphism~$\wn_i$, so that we obtain a formula analogous
to~\eqref{hatHH'H''},
\begin{equation}
\label{HH'H''}
 H_i=H_i^{\p\p-1}\circ H_i'\circ H_i'',
\end{equation}
where $H_i''$~is a homeomorphism from~$\Ll_{t_0}'$ onto~$\Ll_{d_i}'$
arising from an isotopy above the path~$\delta_i$ of
Definition~\ref{t:loops}. As above, we shall ask~$H_i''$ to be linked to
the homeomorphism~$\hat H_i''$ of formula~\eqref{hatHH'H''} by a diagram
similar to diagram~\eqref{HPhi}. We notice that $H_i''$~leaves fixed the
points of~$\Mm'$ since the isotopy above~$\delta_i$ giving rise to it
satisfies condition~(III) of Lemma~\ref{t:G}. Following
Notation~\ref{t:notn1}, item~\ref{t:notn1.topalg}, we denote by~$\bar
H_{i\lstar}''$ the isomorphism induced by~$H_i''$ between
$H_n(\Ll_{t_0}',\Mm')$ and $H_n(\Ll_{d_i}',\Mm')$ to distinguish it from the
isomorphism~$H_{i\lstar}''$ induced between
$H_n(\Ll_{t_0}')$ and $H_n(\Ll_{d_i}')$. The link
between~$V_i$ and~$V_i'$ is then given by the next lemma.

\begin{lem}
\label{t:VV'}
The following diagram is commutative:
\begin{equation*}
\begin{CD}
 H_n(\Ll_{d_i}',\Mm')  @>V_i'>>
    H_n(\Ll_{d_i}')\\
 @AA\vis\bar H_{i\lstar}''A
    @AA\vis H_{i\lstar}''A\\
 H_n(\Ll_{t_0}',\Mm')  @>V_i>>
    H_n(\Ll_{t_0}').
\end{CD}
\end{equation*}
\end{lem}

\begin{proof}
This is a straightforward check using the definitions of $V_i$ and~$V_i'$
and formula~\eqref{HH'H''}.
\end{proof}

Next, we make, as earlier, a reduction to the bundle
$P^{\p-1}(\partial\Delta_i)$ for the right-hand side of
diagram~\eqref{Outer}, this time by going back to the definition of
isomorphism~$\wn_i$. We consider the following diagram:
\begin{equation}
\label{Phiwn}
\begin{CD}
 \makebox[0pt]{$H_n(P^{\p-1}(\partial\Delta_i))$}  @.
   \darrow[is]{33}{-15}{-12}{134}{-4}{13}{\wn_i}\\
 @AA\incl_*A\\
 H_n(\Ll_{d_i}^\shrpp)
  @<\shis{\hat H_{i\lstar}''}<<
     H_n(\Ll_{t_0}^\shrpp) @>\can>>
       H_n(\Ll_{t_0}^\shrpp)/\im(\hat H_{i\lstar}-\id)\\
 @VV\vis\Phi_* V   @VV\vis\Phi_* V   @VV\vis\overline{\Phi_*} V\\
 H_n(\Ll_{d_i}')
  @<\shis{H_{i\lstar}''}<<
     H_n(\Ll_{t_0}') @>\can>>
       H_n(\Ll_{t_0}')/\im(H_{i\lstar}-\id).
\end{CD}
\end{equation}
The upper triangle is commutative by the very definition of
isomorphism~$\wn_i$ (cf.~\eqref{Wang} and~\eqref{wn^-1}). The right-hand
square is commutative by the definition of~$\overline{\Phi_*}$ given after
diagram~\eqref{Deg}. Finally, the left-hand square is also commutative,
since we took care of defining~$H_i''$
and~$\hat H_i''$ coherently.

Finally, we come to the left side of diagram~\eqref{Outer}. We have the
following diagram:
\begin{equation}
\label{deldel}
\begin{CD}
 \darrow{46}{-15}{28}{58}{-10}{4}{\partial}
   @.  H_{n-1}(\Mm')\\
 @.  @AA\partial A\\
 H_n(\Ll_{t_0}',\Mm')
   @>\shis{\bar H_{i\lstar}''}>>
     H_n(\Ll_{d_i}',\Mm').
\end{CD}
\end{equation}
It is commutative due to the fact that $H_i''$~leaves fixed the
points of~$\Mm'$ as already pointed out.

Using now the commutativity of the diagram of Lemma~\ref{t:VV'} and of
diagrams~\eqref{Phiwn} and~\eqref{deldel} and taking into account that some
of the arrows, as indicated, are isomorphisms, we see that we only need to
prove that the following diagram commutes once reversed the isomorphism
labeled~$\Phi_*$:
\begin{equation}
\label{kappaV}
\begin{CD}
 H_n(\Mm'\times\partial\Delta_i)  @>\incl_*>>
   H_n(P^{\p-1}(\partial\Delta_i))\\
 @AA\kappa_i A  @AA\incl_*A\\
 H_{n-1}(\Mm')  @.  
   H_n(\Ll_{d_i}^\shrpp)\\
 @AA\partial A  @VV\vis\Phi_*V\\
 H_n(\Ll_{d_i}',\Mm')  @>V_i'>>  H_n(\Ll_{d_i}').
\end{CD}
\end{equation}

To show this, we shall work at the chain level and use a cross-product
defined at this level as in~\cite[Notation~4.5]{London}. Let
$c_i\colon\Ll_{d_i}'\to\Ll_{d_i}^\shrpp$ be the inverse isomorphism of the
one induced by~$\Phi$. It will be convenient to
denote~$c_{i\lbullet}(\Gamma)$ by~$\Gamma^\sharp$ for any singular
chain~$\Gamma$ of~$\Ll_{d_i}'$. Then using Notation~\ref{t:notn1}, items
\ref{t:notn1.Hom} and~\ref{t:notn1.topalg}, the commutativity of
diagram~\eqref{kappaV}, with arrow~$\Phi_*$ reversed, amounts to the
following homology: for any relative $n$\nbd cycle~$\Gamma$ on~$\Ll_{d_i}'$
modulo~$\Mm'$,
\begin{equation}
\label{Homology}
 (-1)^{n-1}\partial\Gamma\times\omega_i\sim
 (H_{i\lbullet}'(\Gamma)-\Gamma)^\sharp\quad
 \text{in $P^{\p-1}(\partial\Delta_i)$}.
\end{equation}
To prove this homology, first observe that
\begin{equation*}
 (H_{i\lbullet}'(\Gamma)-\Gamma)^\sharp=
 H_{i\lbullet}'(\Gamma)^\sharp-\Gamma^\sharp=
 \hat H_{i\lbullet}'(\Gamma^\sharp)-\Gamma^\sharp
\end{equation*}
since we have ensured that $H_i'$ and~$\hat H_i'$ commute with the blowing
down morphism. Homology~\eqref{Homology} will then be given by the
isotopy~$\hat G_i'$ giving rise to~$\hat H_i'$. More precisely, if
$\iota$~is the $1$\nbd simplex of~$\I$ consisting of the identity map, then
$\Gamma^\sharp\times\iota$ is a chain of $\Ll_{d_i}^\shrpp\times\I$ to
which we can apply~$\hat G_i'$, obtaining a chain of~$\hat W'$, in fact
of~$P^{\p-1}(\partial\Delta_i)$ by condition~(\^I\^I) of Lemma~\ref{t:G}
(where $t_0$~must be replaced by~$d_i$ and $\gamma_i$~by~$\omega_i$). We
shall show that
\begin{equation}
\label{delG}
 \partial\hat G_{i\lbullet}'(\Gamma^\sharp\times\iota)=
 \partial\Gamma\times\omega_i-
  (-1)^{n-1}(\hat H_{i\lbullet}'(\Gamma^\sharp)-\Gamma^\sharp).
\end{equation}

Here is the computation; it can already be found in~\cite[p.~540]{London} in
the course of the proof of~\cite[Proposition~4.13]{London} but no result is
stated there which we could refer to. We have
\begin{equation*}
 \partial\hat G_{i\lbullet}'(\Gamma^\sharp\times\iota)=
 \hat G_{i\lbullet}'(\partial\Gamma^\sharp\times\iota)+
 (-1)^n\hat G_{i\lbullet}'(\Gamma^\sharp\times\partial\iota).
\end{equation*}
Concerning the first term of this sum, observe that
$\partial\Gamma^\sharp$~is a chain of $\Mm'\times\set{d_i}$ and that
the restriction of~$\hat G_i'$ to
$(\Mm'\times\set{d_i})\times\I$ coincides with that of
$\Phi\times\omega_i$ by condition~(\^I\^I\^I) of Lemma~\ref{t:G}. Then
\begin{equation*}
\begin{split}
 \hat G_{i\lbullet}'(\partial\Gamma^\sharp\times\iota)
 &=(\Phi\times\omega_i)_\bullet(\partial\Gamma^\sharp\times\iota)=
 \Phi_\bullet(\partial\Gamma^\sharp)\times\omega_{i\lbullet}(\iota)\\
 &=\partial\Gamma\times\omega_i.
\end{split}
\end{equation*}
As to the second term of the sum above, let $\hat0$
and~$\hat1$ be the
$0$\nbd simplices of~$\I$ with respective values~$0$ and~$1$.
Then $\Gamma^\sharp\times\partial\iota=\Gamma^\sharp\times\hat1-
\Gamma^\sharp\times\hat0$, a difference of chains of
$\Ll_{d_i}^\shrpp\times\set{1}$ and
$\Ll_{d_i}^\shrpp\times\set{0}$. But, if $\varpi$~is the
projection of $\Ll_{d_i}^\shrpp\times\I$ onto the
first factor, the restriction of~$\hat G_i'$ to
$\Ll_{d_i}^\shrpp\times\set{0}$ is the same as the
restriction of~$\varpi$ to the same space, by condition~(\^I) of
Lemma~\ref{t:G}. And the restriction of~$\hat G_i'$ to
$\Ll_{d_i}^\shrpp\times\set{1}$ is the same as the
restriction of~$\hat H_i'\circ\varpi$, by the definition of~$\hat H_i'$
in Lemma~\ref{t:H}. Hence
\begin{equation*}
\begin{split}
 \hat G_{i\lbullet}'(\Gamma^\sharp\times\partial\iota)
 &=\hat G_{i\lbullet}'(\Gamma^\sharp\times\hat1)-
 \hat G_{i\lbullet}'(\Gamma^\sharp\times\hat0)\\
 &=\hat H_{i\lbullet}'(\varpi_\bullet(\Gamma^\sharp\times\hat1))-
 \varpi_\bullet(\Gamma^\sharp\times\hat0)\\
 &=\hat H_{i\lbullet}'(\Gamma^\sharp)-\Gamma^\sharp
\end{split}
\end{equation*}
since the chain cross-product we have considered has the
property that, for any spaces $E$ and~$F$, the projection $\varpi\colon
E\times F\to E$ acts as
$\varpi_\bullet(\gamma\times\sigma)=\gamma$ for every chain~$\gamma$ of~$E$
and every $0$\nbd simplex~$\sigma$ of~$F$
(cf.~\cite[Notation~4.5]{London}).

This shows equality~\eqref{delG}, proving homology~\eqref{Homology} and
hence the commutativity of diagram~\eqref{kappaV}. The commutativity of
diagram~\eqref{Outer}, which was reduced to the former, follows.

\paragraph{}

The commutativity of diagrams~\eqref{DegWang} and~\eqref{Outer} implies
that of diagram~\eqref{VarDeg} and hence proves
Proposition~\ref{t:vardeg}.\qed

\section{A generalization of the Zariski-van Kampen theorem to higher
homotopy}

We give here a projective version of the van~Kampen type theorem
of~\cite[Theorem 2.4]{Annals} using the above defined homotopy
variation operators (see section~\ref{s:var}) instead of the
degeneration operators of~\cite{Annals}. In fact we shall give two
proofs of this result. One is based on an affine version of the
Zariski-van Kampen type theorem from~\cite{Annals} and the other
depends on Theorem~5.1 from~\cite{London}.

\begin{theo}
\label{t:highvK}
Let $V$~be a hypersurface in~$\PP^{n+1}$ with $n\geq2$ having only
isolated singularities. Consider a pencil $(L_t)_{t \in\PP^1}$ of
hyperplanes in~$\PP^{n+1}$ with the base locus~$\Mm$ transversal
to~$V$. Denote by $t_1$, \dots,~$t_N$ the collection of those~$t$
for which $L_t\cap V$ has singularities. Let $t_0$~be different
from either of $t_1$, \dots,~$t_N$. Let $\gamma_i$~be a good
collection, in the sense of Definition~\ref{t:loops}, of paths
in~$\PP^1$ based in~$t_0$. Let $e\in\Mm-\Mm\cap V$ be a base point.
Let $\Vv_i$~be the variation operator corresponding to~$\gamma_i$.
Then inclusion induces an isomorphism:
\begin{equation*}
 \pi_n(\PP^{n+1}-V,e)\leftisarrow\pi_n(L_{t_0}-L_{t_0}\cap V,e)\Big/
 \sum_{i=1}^N\im\Vv_i.
\end{equation*}
\end{theo}

\paragraph{First Proof.} We apply Theorem~5.1 of~\cite{London} to 
the non-singular quasi-pro\-jective variety $W-j(V)$ in~$\PP^{n+2}$
(cf.~section~\ref{s:ramcov}). The base locus~$\Mm$ of the pencil
$(\Ll_t)_{t\in\PP^1}$ is transversal to the Whitney
stratification~$\Sigma$ of~$W$ adapted to~$j(V)$
(cf.~Claim~\ref{t:Transv}). Hence \cite{London}~gives the following
isomorphism induced by  inclusion: 
\begin{equation*}
 H_n(W-j(V))\leftisarrow
 H_n\bigl(\Ll_{t_0}\cap (W-j(V))\bigr)\Big/\sum_{i=1}^N\im V_i
\end{equation*}
where the~$V_i$ are the homological variation operators defined in
section~\ref{s:var}.

Recall (Lemma~\ref{t:homHom}) that we have an isomorphism~$\eta$:
\begin{equation*}
 H_n(W-j(V))\leftisarrow\pi_n(\PP^{n+1}-V,e).
\end{equation*}
Now the result follows using the isomorphism~$\alpha_{t_0}$ and the
commutative diagrams of Lemma~\ref{t:homHom}, and the definition
of~$\Vv_i$ by means of~$V_i$ from section~\ref{s:var}.\qed
 
\paragraph{Second Proof.} Let $\CC_{t_0}^n$~denote the affine part
of~$L_{t_0}$ (that is~$L_{t_0}-M$). The group
$\pi_{n}(\CC^n_{t_0}-\CC^n_{t_0}\cap V)$, as
in~\cite[section~1]{Annals}, will be viewed as a module over
$\ZZ[\pi_1(\CC^n_{t_0}-\CC^n_{t_0}\cap V)]=\ZZ[s,s^{-1}]$. We can
use the affine monodromy of $\CC^n_{t_0}-\CC^n_{t_0}\cap V$ which is
a restriction of the projective one. For $1\leq i\leq N$, there is
a commutative diagram:

\begin{equation*}
 \begin{CD}
  \pi_{n}(\CC^n_{t_0}-\CC^n_{t_0}\cap V)
  @>h_i>>
  \pi_{n}(\CC^n_{t_0}-\CC^n_{t_0}\cap V)\\
  @VVs^d-1V   @VVs^d-1V\\
  \pi_{n}(\CC^n_{t_0}-\CC^n_{t_0}\cap V)
  @>h_i>>
  \pi_{n}(\CC^n_{t_0}-\CC^n_{t_0}\cap V)\\
  @VVV   @VVV\\
  \coker(s^d-1)
  @>h_i>>
  \coker(s^d-1)\\
  @VVa_*V   @VVa_*V\\
  \pi_{n}(L_{t_0}-L_{t_0}\cap V)
  @>h_i>>
  \pi_{n}(L_{t_0}-L_{t_0}\cap V),
 \end{CD}
\end{equation*}
where $a_*$~is an isomorphism of $\ZZ[s,s^{-1}]$\nbd modules
(cf.~Lemma~1.13 in~\cite{Annals}). Hence $a_*$~identifies the image
of the affine monodromy with the image of the projective one. A
similar diagram yields the equality of the images of degeneration
operators. Now the theorem follows from
Proposition~\ref{t:vardeg}.\qed

\begin{rem}
\label{generalize}
We presented Theorem~\ref{t:highvK} as a generalization of
the classical Zariski-van~Kampen theorem. However the latter
concerns the case $n=1$ and our definition of the homotopy variation
operators does not make sense in this case. But, as $M$~is then
reduced to a point, Lemma~\ref{t:varmonodr} makes it natural to
consider that in this case ${\cal V}_i$~should be nothing more than
$h_{i\#}-{\rm id}$. This in turn does not make sense since
$\pi_1(L_{t_0}-L_{t_0}\cap V,e)$ is not commutative but amounts to
saying that the only identifications to make are of each~$x$ of
$\pi_1(L_{t_0}-L_{t_0}\cap V,e)$ with each~$h_{i\#}(x)$. Our
theorem then actually reduces to the classical Zariski-van Kampen
theorem. Nevertheless our proof does not work in the case $n=1$. The
statements generalize, but not the proofs.
\end{rem}

\end{document}